\documentclass[a4paper,14pt]{article}

\usepackage{mathrsfs}
\usepackage{amsmath,amscd}
\usepackage{amssymb, amsmath}
\usepackage{graphicx}
\usepackage{color}

\allowdisplaybreaks

\newcommand{\newcom}{\newcommand}

\newcom{\cA}{{\mathcal A}}
\newcom{\cB}{{\mathcal B}}
\newcom{\cC}{{\mathcal C}}
\newcom{\cD}{{\mathcal D}}
\newcom{\cE}{{\mathcal E}}
\newcom{\cF}{{\mathcal F}}

\newcom{\cJ}{{\mathcal J}}

\newcom{\cL}{{\mathcal L}}
\newcom{\cM}{{\mathcal M}}
\newcom{\cP}{{\mathcal P}}
\newcom{\cS}{{\mathcal S}}
\newcom{\cQ}{{\mathcal Q}}
\newcom{\cT}{{\mathcal T}}
\newcom{\cY}{{\mathcal Y}}
\newcom{\cZ}{{\mathcal Z}}
\newcom{\R}{\mathbb R}
\newcom{\T}{\mathbb T}
\newcom{\N}{\mathbb N}
\newcom{\Z}{\mathbb Z}
\newcom{\C}{\mathbb C}
\newcom{\E}{\mathbb E}

\newcom{\e}{\epsilon}

\newcom{\al}{\alpha}
\newcom{\be}{\beta}
\newcom{\del}{\delta}

\newcom{\ga}{\gamma}
\newcom{\Ga}{\Gamma}

\newcom{\Lam}{\Lambda}
\newcom{\lam}{\lambda}

\newcom{\Om}{\Omega}
\newcom{\om}{\omega}

\newcom{\Si}{\Sigma}
\newcom{\si}{\sigma}
\newcom{\s}{\varsigma}

\newcom{\tht}{\theta}
\newcom{\dtri}{\nabla}
\newcom{\tri}{\triangle}

\newcom{\f}{\frac}
\newcom\na{\nabla}
\newcom{\Del}{\Delta}
\newcom{\ep}{\epsilon}

\newcom{\p}{\partial}

\newcom{\uep}{{\bf u}_{\epsilon}}
\newcom{\nep}{n_{\epsilon}}
\newcom{\cep}{c_{\epsilon}}

\newcom{\beq}{\begin{equation}}
\newcom{\eeq}{\end{equation}}
\newcom{\ben}{\begin{eqnarray}}
\newcom{\een}{\end{eqnarray}}
\newcom{\beno}{\begin{eqnarray*}}
\newcom{\eeno}{\end{eqnarray*}}
\newcom{\bal}{\begin{aligned}}
\newcom{\eal}{\end{aligned}}

\numberwithin{equation}{section}

\title{Small Solitons and Multi-Solitons in Generalized
Davey-Stewartson System}

\author{Mengxue Bai$^{1}$\thanks{E-mail: mengxuebai@163.com},
\quad  Jian Zhang$^{1}$\thanks{Corresponding author. E-mail: zhangjian@uestc.edu.cn}~,
\quad  Shihui Zhu$^{2}$\thanks{E-mail: shihuizhumath@163.com}\\
\emph{\small  $^{1}$School of Mathematical Sciences, University of Electronic Science and Technology of China, }\\
\emph{\small Sichuan, Chengdu 611731, China}\\
\emph{\small  $^{2}$School of Mathematical Sciences, Sichuan Normal University, }\\
\emph{\small Sichuan, Chengdu 610066, China}}

\date{}

\begin{document}
\maketitle
{\textbf{Abstract:}  This paper is concerned with the generalized Davey-Stewarston system in two dimensional space. Existence and stability of small solitons are proved by solving two correlative constrained variational problems and spectrum analysis. In addition, multi-solitons with different speeds are constructed by bootstrap argument.}



\section{Introduction}
\renewcommand{\theequation}{\thesection.\arabic{equation}}

\par   Consider the generalized Davey-Stewartson system in two dimensional space,
\begin{align}
i\varphi_{t}+\Delta \varphi+\lvert\varphi\rvert^{p-1}\varphi+E_{1}(\lvert\varphi\rvert^{2})\varphi=0,\;\;(t,\;x)\in\mathbb{R}\times\mathbb{R}^{2}.
\end{align}
Here $1 < p <\infty$ and $E_{1}$ is the singular integral operator with symbol $\sigma_{1}(\xi)=\frac{\xi^{2}_{1}}{\lvert\xi\rvert^{2}}$, $\xi\in\mathbb{R}^{2}$, that is $E_{1}(\varphi)(x)=\mathcal{F}^{-1}(\frac{\xi^{2}_{1}}{\lvert\xi\rvert^{2}}\mathcal{F}(\varphi)(\xi))$, where $\mathcal{F}$ and $\mathcal{F}^{-1}$ represent the Fourier transform and Fourier inverse transform on $\mathbb{R}^{2}$ respectively, and $\mathcal{F}(\varphi)(\xi)=\frac{1}{2\pi}\int e^{-ix\xi}\varphi(x)dx$. Here and hereafter we denote $\int_{\mathbb{R}^{2}}\cdot dx$ by $\int\cdot dx$, $L^{2}(\mathbb{R}^{2})$ by $L^{2}$ and $H^{1}(\mathbb{R}^{2})$ by $H^{1}$.
\par  (1.1) origins from fluid mechanics, and it models the evolution of weakly nonlinear water waves having a predominant direction of travel. More precisely, (1.1) is the extension of the Davey-Stewartson systems in the elliptic-elliptic case, namely
\begin{equation}
\left\{\begin{array}{lll}
& i\varphi_{t}+
\lambda\varphi_{x_1x_1}+\mu
\varphi_{x_2x_2}=a\lvert\varphi\rvert^{2}\varphi+b_{1}\varphi \phi_{x_1},\\
 & \nu\phi_{x_1x_1}+\phi_{x_2x_2}=-b_{2}(\lvert\varphi\rvert^2)_{x_1}
\end{array}\right.
\end{equation}
\\
($a\in \mathbb{R},\;\lambda,\;\mu,\;\nu,\;b_{1}\;and\;b_{2} > 0$) which describes the time evolution of two-dimensional surface of water wave having a propagation preponderantly in the $x_1$-direction (see \cite{C1992, GS1990, GW1999, O1994}).
\par  Ghidaglia and Saut\cite{GS1990} showed the local well-posedness of the Cauchy problem of (1.1) in the natural energy space $H^{1}$ for $p=3$, then Guo and Wang \cite{GW1999} generalized this result to $1 < p < \infty$. Ozawa \cite{O1992} constructed the exact blow up solutions of the Cauchy problem of (1.1) for $ p = 3$  (also see the numerical simulation result of Sulem C. and Sulem P. L.\cite{SS1999}). By Ghidaglia and Saut\cite{GS1990} as well as  Ohta \cite{O1995B}, it was known that the Cauchy problem of (1.1) has blow up solutions to appear for $1 < p < \infty$. In addition, Gan and Zhang \cite{GZ2008} studied sharp threshold of blow up and global existence for the Cauchy problem of (1.1). In terms of Zhang's argument\cite{Z2000}, Zhu \cite{Z2016} got global existence of small solutions with the mass for the Cauchy problem of (1.1).
\par For $\omega > 0$, consider the following nonlinear elliptic equation
\begin{align}
\Delta u+\lvert u\rvert^{p-1}u+E_{1}(\lvert u\rvert^{2})u=\omega u,\;\;u\in H^{1}
\end{align}
If $u(x)$ is a non-trivial solution of (1.3), then $e^{i\omega t}u(x)$ is a soliton of (1.1).
\par  Cipolatti \cite{C1992} proved the existence of positive solutions of (1.3) by means of P. L. Lion's concentration-compactness method (see \cite{L1984A, L1984B}). Then Cipolatti (see \cite{C1993}), Ohta \cite{O1995B}, Gan and Zhang \cite{GZ2008} showed the instability of the solitons of  (1.1) for $3 \leq p < \infty$ respectively by different methods.  Because of the singular operator $E_{1}$ in (1.3) (see \cite{SS1999}), the uniqueness of positive solutions for (1.3) is still open. Under the assumption of uniqueness of positive solutions of (1.3), Ohta \cite{O1994} proved that for $1 < p < 3$,  there exists a sequence of frequency  $\omega_{n} > 0$ such that $\omega_{n} \rightarrow 0$ and the solitons $e^{i\omega_n t}Q_{n}$ are stable, where $Q_{n}$ is the unique positive solution of (1.3) corresponding to $\omega_{n} > 0$. Moreover Ohta \cite {O1995B} got the stability of the solitons generated by the set of minimizers of the associated variational problem. From \cite{C1993, O1995A, GZ2008}, the instability of solitons for (1.1) has gotten a comprehensive study. And from \cite{O1994, O1995B}, further study to stability of solitons for (1.1) becomes an interesting topic. In this paper we develop some new technologyies to study stability of solitons for (1.1).
\par For $u\in H^{1}\backslash\{0\}$, we define the funtional
\begin{align}
J(u)=\frac{(\int\lvert u\rvert^{2}dx)(\int\lvert \nabla u\rvert^{2}dx)}{\int E_{1}(\lvert u\rvert^{2})\lvert u\rvert^{2}dx}.
\end{align}
Then we  consider the variational problem
\begin{align}
d_{J}=inf_{\{u\in H^{1}\backslash\{0\}\}}J(u).
\end{align}
 It is known that (1.5) possesses a positive minimizer $u\in H^{1}$ (see \cite{Z2016}). Therefore $d_{J}$ is a positive constant. Moreover, for arbitrary $u\in H^{1}$, one has the sharp interpolation inequality:
\begin{align}
\int E_{1}(\lvert u\rvert ^2)\lvert u\rvert^2dx\leq
\frac{1}{d_{J}}\int\rvert\nabla u\lvert^{2}dx\int\rvert u\lvert^2dx.
\end{align}
Let
%\begin{align}
%H_{J}=\{\phi\in H^{1}(\mathbb{R}^{d}),\;\;\int\lvert\phi\rvert^{2}dx<2d_{J}\},
%\end{align}
%and
\begin{align}
\omega_{J}=sup\{\omega\in \mathbb{R}\big\lvert\;\|Q_{\omega}\|_{L^{2}} < \sqrt{2d_{J}}\},
\end{align}
where $Q_{\omega}$ is the positive solution of (1.3). Firstly we can prove that $\omega_{J} > 0$. Then we prove the following crucial results of stability of solitons for (1.1).
\par \noindent\textbf{Theorem A.} Let $\omega\in (0, w_{J})$, $1 < p < 3$ and $Q_{\omega}(x)$ is the positive solution of (1.3). Suppose that the positive solution of (1.3) is unique, then the small solitons $e^{i\omega t}Q_{\omega}(x)$ of (1.1) is orbitally stable. Moreover it is true that $\frac{d}{d\omega}\int Q^{2}_{\omega}dx > 0$ for all $\omega\in (0, \omega_{J})$.
\par  In order to prove Theorem A, we construct and solve two correlative constrained variational problems. Then we ascertain frequency from mass by establishing a one-to-one mapping. Finally we bridge Grillakis-Shatah-Strauss method \cite{GSS1987} and Cazenave-Lions method  \cite{CL1982} for stability of the solitons  by spectrum analysis. It is clear that Theorem A includes the results in \cite{O1994, O1995B}. Moreover
technologies developed in this paper can be used to determine frequency from the prescribed mass in the  normalized solution problems (see \cite{BJS2016, BMRV2021, BZZ2020, Z2000}). We discuss this problem in other papers.
\par
In terms of C$\widehat{o}$te and Le Coz's arguments \cite{CC2011}, Wang and Cui \cite{WC2017} constructed the high speed excited multi-solitons of (1.1). Multi-solitons  are concerned with the famous soliton resolution conjecture, which is emphasized in  Tao \cite{T2009}, Zakharov and Shabat \cite{Z1972}. The stability of solitons and the soliton resolution problems are crucial topics in understanding the dynamics of nonlinear dispersive evolution equations (see Tao \cite{T2006}). Therefore we use the stable solitons of (1.1) obtained in Theorem A to construct multi-solitons with different speeds for $(1.1)$ according to Martel, Merle and Tsai's scheme (see \cite{MM2006, MMT2006}). %It is well known that C$\widehat{o}$te and Le Coz's arguments \cite{CC2011} originate from Martel, Merle and Tsai's scheme (see \cite{MM2006} and \cite{MMT2006}).
We prove the following theorem.
\par \noindent\textbf{Theorem B: } Let $1 < p < 3$. For $K \geq 2$ and $k=1,\;2,\;\cdot\cdot\cdot,\;K$, taking $\omega_{k}\in (0,\;\omega_{J})$, $\gamma_{k}\in \mathbb{R}$, $x_{k}\in \mathbb{R}^{2}$, $v_{k}\in \mathbb{R}^{2}$ with $v_{k}\neq v_{k'}$ to $k\neq k'$ and
\begin{eqnarray}
R_{k}(t,x)=Q_{\omega_{k}}(x-x_{k}-v_{k}t)e^{i(\frac{1}{2}v_{k}x-\frac{1}{4}\lvert v_{k}\rvert^{2}t+\omega_{k}t+\gamma_{k})}
\end{eqnarray}
with $(t,x)\in\ \mathbb{R}\times \mathbb{R}^{2}$, there exists a solution $\varphi(t,\;x)$ of (1.1) such that
\begin{align}
\lim\limits_{t \rightarrow +\infty}\|\varphi(t)-\sum_{k=1}^{K}R_{k}(t)\|_{H^{1}}=0.
\end{align}
The solution $\varphi(t,\;x)$ of (1.1) holding (1.9) is called multi-soliton of (1.1).
\par  The rest contents of this paper are organized as follows. In section 2, we show global existence of small solutions  of the Cauchy problem for (1.1) and existence of solitons for (1.1). In section 3, by solving two correlative constrained variational problems, we establish a one-to-one mapping between mass and frequency. In section 4,
we prove orbital stability of small solitons depending on mass for (1.1). Moreover by spectrum analysis, we communicate the relationship between Cazenave-Lions method \cite{CL1982} and Grillakis-Shatah-Strauss method \cite{GSS1987}. In addition, we get orbital stability of solitons depending on  frequencies $\omega\in (0, \omega_{J})$ for (1.1). In section 5, we construct multi-solitons with different speeds for (1.1) by all stable solitons in terms of the bootstrap scheme and the uniform backward estimate according to \cite{CC2011, MM2006, MMT2006, WC2017}.
\section{Well-Posedness}
For $t_{0}\in \mathbb{R}$, we impose the initial data of (1.1) as follows.
\begin{align}
\varphi(t_{0},x)=\varphi_{0}(x),\;\;\;x\in \mathbb{R}^{2}.
\end{align}
In $H^{1}$, we define the energy functional
\begin{align}
E(\varphi)=\int\lvert \nabla \varphi\rvert^{2}dx-\frac{2}{p+1}\int\lvert \varphi\rvert^{p+1}dx-\frac{1}{2}\int E_{1}(\lvert \varphi\rvert^{2})\lvert \varphi\rvert^{2}dx;
\end{align}
the mass functional
\begin{align}
M(\varphi)=\int\lvert \varphi\rvert^{2}dx;
\end{align}
and the momentum functional
\begin{align}
P(\varphi)=Im\int \overline{\varphi}\nabla \varphi dx.
\end{align}
\par  First we have the following lemma by Zhu \cite{Z2016}.
\par \noindent\textbf{Lemma 2.1.} Define the variational problem
\begin{align}
d_{J}=inf_{\{\phi\in H^{1}\backslash\{0\}\}}J(\phi)\;\;\;with\;J(\phi)=\frac{(\int\lvert \phi \rvert^{2}dx)(\int\lvert \nabla\phi \rvert^{2}dx)}{\int E_{1}(\lvert \phi \rvert^{2})\lvert \phi \rvert^{2}dx}.
\end{align}
Then  (2.5) possesses a nontrivial minimizer and $d_{J} > 0$. Moreover for $\phi \in H^{1}$,
we have the sharp Gagliardo-Nirenberg type inequality
\begin{align}
\int E_{1}(\lvert \phi \rvert^{2})\lvert \phi \rvert^{2}dx \leq \frac{1}{d_{J}}\int\rvert\nabla\phi\lvert^{2}dx\int\rvert\phi\lvert^{2}dx.
\end{align}
%Now we define the space,
%\begin{align}
%H_{J}=\{\phi\in H^{1}(\mathbb{R}^{d}),\;\;\int\lvert\phi\rvert^{2}dx<2d_{J}\}.
%\end{align}
Then we have the following global well-posedness  for (1.1) with small mass.
\par \noindent\textbf{Theorem 2.2.} Let $1 < p < 3$, $\varphi_{0}\in H^{1}$ and $\|\varphi_{0}\|_{L^{2}} < \sqrt{2d_{J}}$. Then the Cauchy problem (1.1)-(2.1) possesses a unique global solution $\varphi(t,x)\in C(\mathbb{R},\;H^{1})$ with mass conservation $M(\varphi)=M(\varphi_{0})$, energy conservation $E(\varphi)=E(\varphi_{0})$ and momentum conservation $P(\varphi)=P(\varphi_{0})$ for all $t\in \mathbb{R}$.
\\
\noindent\textbf{Proof.} By \cite{C2003} and \cite{GV1979}, for $\varphi_{0}\in H^{1}$ with $\|\varphi_{0}\|_{L^{2}} < \sqrt{2d_{J}}$, there exists a unique solution $\varphi (t,\;x)$ of the Cauchy problem (1.1)-(2.1) in $C((-T,\;T);H^{1})$ to some $T>0$ (maximal existence time). And $\varphi(t,\;\cdot)$ satisfies mass conservation $M(\varphi)=M(\varphi_{0})$, energy conservation $E(\varphi)=E(\varphi_{0})$ and momentum conservation $P(\varphi)=P(\varphi_{0})$ for all $t\in (-T,\;T)$. Furthermore one has the alternatives: $T=\infty$ (global existence) or else $T<\infty$ and $lim_{t\rightarrow T}\|\varphi\|_{H^{1}}=\infty$ (blow up). Thus from (2.2), (2.6), we have that
\begin{align}
E(\varphi) \geq (1-\frac{1}{2d_{J}}\int \lvert \varphi \rvert^{2}dx)\int \lvert \nabla\varphi \rvert^{2}dx-\frac{2}{p+1}\int \lvert \varphi \rvert^{p+1}dx.
\end{align}
From the Gagliardo-Nirenberg inequality
\begin{align}
\int\lvert \varphi \rvert^{p+1}dx \leq C(p)(\int\lvert \varphi \rvert^{2}dx)(\int\lvert \nabla\varphi \rvert^{2}dx)^{\frac{p-1}{2}},\;\;\varphi\in H^{1},
\end{align}
mass conservation and energy conservation, (2.7) yields that
\begin{align}
C_{1}\int\lvert\nabla \varphi\rvert^{2}dx-C_{2}(\int\lvert\nabla \varphi\rvert^{2}dx)^{\frac{p-1}{2}}\leq E(\varphi_{0}),
\end{align}
where $C_{1}$ and $C_{2}$ are positive constants only concerning $d$ and $\varphi_{0}\in  H^{1}$. From $1 < p < 3$, $\int\lvert\nabla \varphi\rvert^{2}dx$ is bounded for $t\in(-T,\;T)$ with any $T<\infty$. Therefore combining with the mass conservation, we get that $\varphi(t,\;x)$ globally exists in $t\in (-\infty,\;\infty)$. Moreover, the mass conservation and the energy conservation, as well as the momentum conservation are true to all $t\in \mathbb{R}$.
\par  This proves Theorem 2.2.
\par  \noindent\textbf{Theorem 2.3.}  Let $1<p<3$ and $\omega\in\mathbb{R}$. Then the necessary condition for the nonlinear elliptic equation
\begin{align}
\Delta u-\omega u+\lvert u\rvert^{p-1}u+E_{1}(\lvert u\rvert^{2})u=0,\;\;u\in H^{1}
\end{align}
to possess nontrivial solutions is $\omega > 0$.
\par  \noindent\textbf{Proof.} Let $u(x)$ be a nontrivial solution of (2.10). By the Pohozaev's identity (see \cite{P1965}), we have that,
\begin{align}
-\frac{2}{p+1}\int \rvert u \lvert^{p+1}dx-\frac{1}{2}\int E_{1}(\rvert u \lvert^{2})\rvert u \lvert^{2}dx+ \omega  \int \rvert u \lvert^{2}dx=0.                \end{align}
Since
\begin{align}
\int E_{1}(\rvert u \lvert^{2})\rvert u \lvert^{2}dx=\int \rvert u \lvert^{2}\mathcal{F}^{-1}(\sigma_{1}(\xi)\mathcal{F}(\rvert u \lvert^{2}))dx=\int \sigma_{1}(\xi)\rvert \mathcal{F}(\rvert u \lvert^{2})\lvert^{2} d\xi > 0,
\end{align}
from (2.11) it follows that $\omega > 0$.
\par  This proves Theorem 2.3.
\par  From Cipolatti \cite{C1992} we state the following lemma.
\par \noindent\textbf{Lemma 2.4.}  Let $1 < p < 3$ and $\omega > 0$. Then the nonlinear elliptic equation (2.10) possesses a positive solution $Q_{\omega}(x)$, and $Q_{\omega}(x)$ has exponential decay property with $C_{1}, C_{2} > 0$:
\begin{align}
\lvert\nabla Q_{\omega}(x)\rvert+\lvert Q_{\omega}(x)\rvert\leq C_{1}e^{-C_{2}\lvert x\rvert},\;\;\;x\;\in\;\mathbb{R}^{2}.
\end{align}
In addition, let $E_{j},\;j=1,\;2$ be the pseudo-differential operator with symbol $\sigma_{j}(\xi)=\frac{\xi_{1}\xi_{j}}{\lvert\xi\rvert^{2}}$. Then $E_{j}(\lvert Q_{\omega}(x)\rvert^{2})$ has exponential decay property:
\begin{align}
\lvert E_{j}(\lvert Q_{\omega}(x)\rvert^{2})\rvert \leq C_{1}e^{-C_{2}\lvert x\rvert},\;\;\;x\;\in\;\mathbb{R}^{2}.
\end{align}
\noindent\textbf{Proof.} The proof of (2.13) is from \cite{C1992}. For reader's convenience, we give the proof of (2.14) (also see \cite{WC2017}). Let $f=\mathcal{B}\ast \lvert Q_{\omega}(x)\rvert^{2}$, where $\mathcal{B}$ is the fundamental solution of the Laplacian. Then $f$ is a solution of the following equation:
\begin{align}
-\Delta f=\lvert Q_{\omega}(x)\rvert^{2}.
\end{align}
It is easy to see that
\begin{align}
E_{j}(\lvert Q_{\omega}(x)\rvert^{2})=-\partial_{1}\partial_{j}f,\;\;j=1,\;2.
\end{align}
Hence, in order to prove (2.14), it is sufficient to prove the spatial exponential decay of $\partial_{1}\partial_{j}f$. Note that
\begin{align*}
-\Delta\partial_{1}\partial_{j}f=2Re(\partial_{1}\partial_{j}Q_{\omega}\overline{Q}_{\omega}+\partial_{1}Q_{\omega}\partial_{j}\overline{Q}_{\omega}),\;\;j=1,\;2.
\end{align*}
By (2.13), we see that there exists positive constants $C_{1}$ such that the absolute value of the right-hand side of the above equation is bounded by $C_{1}e^{-C_{2}\lvert x\rvert}$, let $g=g(\lvert x\rvert)$ be the unique radial solution of the problem
\begin{align*}
-\Delta g=C_{1}e^{-\sqrt{\omega}\lvert x\rvert},\;\;\;\;\lim_{\lvert x\rvert\rightarrow \infty} g(\lvert x\rvert)=0.
\end{align*}
A simple computation shows that there exists a polynomial $P(x)$ such that
\begin{align*}
\lvert g(\lvert x\rvert)\rvert \leq C_{1}P(x)e^{-C_{2}\lvert x\rvert},\;\;for\;\lvert x\rvert \geq 0.
\end{align*}
Hence, by using the standard super and sub-solutions method, we obtain $\lvert\partial_{1}\partial_{j}f\rvert \leq C_{1}e^{-C_{2}\lvert x\rvert}$.
\par  This completes the proof.

\section{Correlative  Variational Framework}
Firstly we state the profile decomposition theory of a bounded sequence in $H^1$, which is proposed by Hmidi and Keraani in \cite{TS2005}.
\par \noindent\textbf{Lemma 3.1.} \label{profile} Let
$\{u_n\}_{n=1}^{+\infty}$ be a bounded sequence in $H^1$. Then,
 there exists a subsequence of $\{u_n\}_{n=1}^{+\infty}$ (still
denoted by $\{u_n\}_{n=1}^{+\infty}$ ) and a sequence
$\{U^j\}_{j=1}^{+\infty}$ in $H^1$ and a family of
$\{x_n^j\}_{j=1}^{+\infty}\subset \mathbb{R}^2$ satisfying the following.
\\
(i) For every $j\neq k$, $\lvert x_n^j-x_n^k\rvert\rightarrow+\infty$  as $n\rightarrow+\infty$.
\\
(ii) For every $l\geq 1$ and every $x\in \mathbb{R}^2$, $u_n(x)$ can be decomposed by
\begin{align*}
 u_n(x)=\sum\limits_{j=1}^{l}U^j(x-x_n^j)+u_n^l,
 \end{align*}
with the remaining term $u_n^l:=u_n^l(x)$ satisfying
\begin{align*}
\lim_{l\rightarrow+\infty}\limsup_{n\rightarrow+\infty}\;\|u_n^l\|_{L^{q}(\mathbb{R}^2)}=0,\;for\;every\;q\in(2,+\infty).
\end{align*}
Moreover, as $n\rightarrow+\infty$,
\begin{align}
\|u_n\|_{L^{2}}^2=\sum\limits_{j=1}^{l}\|U^j\|_{L^{2}}^2+\|u_n^l\|_{L^{2}}^2+o(1),\; \|\nabla
u_n\|_{L^{2}}^2=\sum\limits_{j=1}^{l}\|\nabla U^j\|_{L^{2}}^2+\|\nabla
u_n^l\|_{L^{2}}^2+o(1),
\end{align}
where $\lim\limits_{n\rightarrow+\infty}o(1)=0$.
\par  The sequence $\{x_n^j\}_{n=1}^{+\infty}$ is called to satisfy the orthogonality condition if and only if for every $k\neq j$,
$\lvert x_n^k-x_n^j\rvert\rightarrow+\infty\;\;as\;\;n\rightarrow \infty$.
\par  Then we show the following lemma.
\par \noindent\textbf{Lemma 3.2.} Let $\{U^{j}\}_{j=1}^{l}$ be a family of bounded sequences in $H^1$ and $\{x_n^j\}_{n=1}^{+\infty}$ be a  orthogonality sequence in $\mathbb{R}^2$. We claim that for every $1 < p < \infty$,
\begin{align}
\int(\sum_{j=1}^{l}U^{j}(x-x_n^j))^{p+1}dx\rightarrow\sum_{j=1}^{l}\int (U^{j}(x-x_n^j))^{p+1}dx\
\ \ { \rm as}\ \ \ n\rightarrow +\infty.
\end{align}
\begin{align}
\int
 E_{1}(\lvert\sum\limits_{j=1}^l U_n^j\rvert^2)  \lvert\sum\limits_{j=1}^l U_n^j\rvert^2dx\rightarrow \sum\limits_{j=1}^l\int
 E_{1}(\lvert U_n^j\rvert^2)  \lvert U_n^j\rvert^2dx\
\ \ { \rm as}\ \ \ n\rightarrow +\infty.
\end{align}
\noindent\textbf{Proof.}  We give the proof of (3.2). Then (3.3) can be obtained by the same arguments, (also see \cite{Z2016}).
Assume that every $U^j$ is continuous and compactly supported. From the basic inequality: for every $p > 1$
\[\lvert \lvert\sum\limits_{j=1}^la_j\rvert^{p+1}-\sum\limits_{j=1}^l\lvert a_j\rvert^{p+1} \rvert\leq
C\sum\limits_{j\neq k}\lvert a_j\rvert \lvert a_k\rvert^p,\] we have that it is sufficient to  prove  that  the mixed terms in the left hand side of  (3.2) vanish.
More precisely, for all $j\neq k$, we claim that
\begin{align}
\int \lvert U^j\rvert \lvert U^k\rvert\lvert U^m\rvert^{p-1}dx\rightarrow 0\;\;as\;\;
n\rightarrow+\infty.
\end{align}
To show (3.4), based on some basic computations we deduce the following inequality
\begin{align}
 \int \lvert U^jU^k\rvert\lvert U^m\rvert dx&\leq  C \left(\int
\lvert U^jU^k\rvert^{\frac{p+1}{2}}dx\right)^{\frac{2}{p+1}}\int\rvert\nabla
U^m\lvert^{p-1}dx\notag\\
&\leq C\left(\int
\lvert U^jU^k\rvert^{\frac{p+1}{2}}dx\right)^{\frac{2}{p+1}}.
\end{align}
 From  Lemma 3.1,
we deduce that
\begin{align}
\int
\lvert U^jU^k\rvert^{\frac{p+1}{2}}dx=\int \lvert U^j(y-(x_n^j-x_n^k))U^k(y)\rvert^{\frac{p+1}{2}}dy\rightarrow 0\ \ {\rm as }\  \ n\rightarrow+\infty.
\end{align}
Finally, from  (3.5) and (3.6), (3.4) can be obtained.
\par  This completes the proof of Lemma 3.2.
\par \noindent\textbf{Theorem 3.3.} Let $1 < p < 3$ and $0<m<2d_{J}$, where $d_{J}$ is defined as Lemma 2.1. We set the constrained variational problem
\begin{align}
d_{m}:=inf_{\{u\in H^{1},\;\int \lvert u\rvert^{2}dx=m\}}E(u).
\end{align}
Then (3.7) possesses a positive minimizer $Q_{m} \in H^{1}$. Moreover there exists a unique $\omega_{m}>0$ such that $Q_{m}$ is the solution of (2.10) with $\omega=\omega_{m}$.
\par \noindent\textbf{Proof.}  It is obvious that $\{u\in H^{1},\;\int \lvert u\rvert^{2}dx=m\}$ is not empty. In the following we complete this proof by four steps.
\par  Step 1. $-\infty < d_{m} < 0 $.
\par  From (2.6), (2.7), (2.8) and the Young inequality, we deduce that
\begin{align}
 E(u)&=\int\lvert\nabla u\rvert^{2}dx-\frac12\int
 E_{1}(\lvert u\rvert^2)\lvert u\rvert^2dx-\frac {2}{p+1} \int\lvert u\rvert^{p+1}dx\notag\\
&\geq (1-\frac{\|u\|_{L^{2}}^2}{2d_{J}})\ \|\nabla u\|_{L^{2}}^2-C\|u\|_2^2\ \|\nabla u\|^{p-1}_{L^{2}}\notag\\
&\geq ((1-\frac{\|u\|_{L^{2}}^2}{2d_{J}})-\varepsilon)\|\nabla u\|_{L^{2}}^2-C(\varepsilon,\|u\|_{L^{2}}),
\end{align}
Taking $0<\varepsilon<1-\frac{\|u\|_{L^{2}}^2}{2d_{J}}$, since $0 < \int \rvert u\lvert^{2}dx=m < 2d_{J}$, by (3.8),
we have that
\begin{align}
E(u) \geq -C(\varepsilon,m)=constant>-\infty.
\end{align}
Let $u_{\lambda}=\lambda u(\lambda x)$. We see that $\|u_{\lambda}\|_{L^{2}}^2=\|u \|_{L^{2}}^2=m<2d_{J}$ and
\begin{align}
E( u_{\lambda} )=\lambda^2(\int\lvert\nabla u\rvert^{2}dx-\frac12\int
 E_{1}(\lvert u\rvert^2)\lvert u\rvert^2dx)-\frac {2\lambda^{p-1}}{p+1} \int\lvert u\rvert^{p+1}dx.
\end{align}
From (2.6), it can be obtained that if $\|u \|_2^2=m<2d_{J}$, then
\begin{align}
\int\lvert \nabla u\rvert^{2}dx-\frac12\int
 E_{1}(\lvert u\rvert^2)\lvert u\rvert^2dx\geq C_1 >0.
\end{align}
Moreover, since $1 < p < 3$, there exists a sufficiently small $0 < \lambda <<1$  such that $E( u^{\lambda} )<0$. It follows that
$d_{m} < 0$. Combining with (3.9), we get that $-\infty < d_{m} < 0$
\par  Step 2. Minimizing sequence is bounded in $H^{1}$.
\par  Let $\{u_n\}^{+\infty}_{n=1}$ be a  minimizing sequence of (3.7). Then we have that
\begin{align}
E(u_n)\rightarrow d_m\;\;as\;\; n\rightarrow+\infty,
\end{align}
\begin{align}
\|u_n\|_{L^{2}}^2= m,\;\;n=1,2,\cdot\cdot\cdot.
\end{align}
By (3.12), one has that
\begin{align}
 E(u_n)<d_m+1\;\;as\;\;n\rightarrow+\infty.
\end{align}
Thus, it can be deduced that for all
$0<\varepsilon<1-\frac{ m }{2d_{J}}$,
\[(1-\frac{m}{2d_{J}}-\varepsilon)\|\nabla u_n\|_{L^{2}}^2\leq d_m+1+C(\varepsilon, m).\]
Combining with $0 < \int \rvert u_{n}(x)\lvert^{2}dx < 2d_{J}$, we deduce that $\{u_{n}\}$ is bounded in $H^{1}$. Moreover since $d_{m} < 0$, one can choose a $0 < \delta < -d_{m}$ to satisfy
\begin{align*}
\frac12\int
 E_{1}(\lvert u_n\rvert^2)\lvert u_n\rvert^2dx+\frac {2}{p+1} \int\lvert u_n\rvert^{p+1}dx=\int\lvert \nabla u_n\rvert^{2}dx-E(u_n)\geq -d_m-\delta,
\end{align*}
for $n$ large enough, which implies that
\begin{align}
\frac12\int
 E_{1}(\lvert u_n\rvert^2)\lvert u_n\rvert^2dx+\frac {2}{p+1} \int\lvert u_n\rvert^{p+1}dx\geq C_0.
\end{align}

\par  Step 3. Existence of minimizer.
\par  We apply Lemma 3.1 to the minimizing sequence $\{u_n\}^{+\infty}_{n=1}$. Then there exists a subsequence still denoted by $\{u_n\}^{+\infty}_{n=1}$ such that
\begin{align}
u_n(x)=\sum\limits_{j=1}^lU_n^j(x)+u_n^l,
\end{align}
where $U_n^j(x):=U^j(x-x_n^j)$ and  $u_n^l:=u_n^l(x)$ satisfies
\begin{align}
\lim\limits_{l\rightarrow+\infty}\limsup\limits_{n\rightarrow+\infty} \|u_n^l\|_{L^{q}(\mathbb{R}^{2})}=0\;\;with\;\;q\in(2,+\infty).
\end{align}
Moreover, by Lemma 3.1 and 3.2, we can get the following estimations as $n\rightarrow +\infty$:
\begin{align}
  \|
u_n\|_{L^{2}}^2=\sum\limits_{j=1}^l\| U_n^j \|_{L^{2}}^2+\|
u_n^l\|_{L^{2}}^2+o(1),
\end{align}
\begin{align}
\|\nabla
u_n\|_{L^{2}}^2=\sum\limits_{j=1}^l\|\nabla U_n^j\|_{L^{2}}^2+\|\nabla
u_n^l\|_{L^{2}}^2+o(1),
\end{align}
\begin{align}
\|
u_n\|_{L^{p+1}}^{p+1}=\sum\limits_{j=1}^l\|
U_n^j\|_{L^{p+1}}^{p+1}+\| u_n^l\|_{L^{p+1}}^{p+1}+o(1),
\end{align}
\begin{align}
\int E_{1}(\lvert u_n\rvert^2) \lvert u_n\rvert^2dx =\sum\limits_{j=1}^l\int
 E_{1}(\lvert U_n^j\rvert^2)  \lvert U_n^j\rvert^2dx+\int  E_{1}(\lvert u_n^l\rvert^2)\lvert
 u_n^l\rvert^2dx +o(1).
\end{align}
From (2.2), (3.16) and (3.18)-(3.21), we have that
\begin{align}
E(u_n)=\sum\limits_{j=1}^{l}E(U_n^j )+E(u_n^l)+o(1)\;\;as\; n\rightarrow+\infty.
\end{align}
For $j=1,2,\cdot\cdot\cdot,l$, let $\widetilde{U}_n^j=\lambda_jU_n^j$ and $\widetilde{u}_n^l=\lambda_n^lu_n^l$,where
\begin{align*}
\lambda_j=\frac{\sqrt{m}}{\|U_n^j \|_{L^{2}}}\geq 1,\;\;\;
\lambda_n^l=\frac{\sqrt{m}}{\|u_n^l\|_{L^{2}}}\geq 1.
\end{align*}
It follows that for $j=1,2,\cdot\cdot\cdot,l$,
\begin{align}
\|\widetilde{U}_n^j \|^{2}_{L^{2}}=\|\widetilde{u}_n^l\|^{2}_{L^{2}}=m.
\end{align}
Moreover, from the convergence of $\sum\limits_{j=1}^{l}\|U_n^j \|_{L^{2}}^2$, one has that
there exists $j_0\geq 1$ such that
\begin{align}
\inf\limits_{j\geq
1} \lambda_j^{p-1}-1 = \ \lambda_{j_0}^{p-1}-1 = (\frac{\sqrt{m}}{\|U^{j_0}\|_{L^{2}}})^{p-1}-1 .
\end{align}
Now, we consider the new energy $E(U_n^j)$ and $E(u_n^l)$.  Then we have that
\begin{align}
E(U_n^j)=\frac{E(\widetilde{U}_n^j)}{\lambda_j^2}+\frac{2(\lambda_j^{p-1}-1)}{p+1}\int\rvert U_n^j \rvert^{p+1}dx+\frac{  \lambda_j^{2}-1 }{2}\int
 E_{1}(\lvert U_n^j\rvert^2)\lvert U_n^j\rvert^2dx,
\end{align}
\begin{align}
E(u_n^l)&=\frac{E(\widetilde{u}_n^l)}{(\lambda_n^l)^{2}}+\frac{2(\lambda_n^l)^{p-1}-1)}{p+1}\int\rvert u_n^l\rvert^{p+1}dx+\frac{(\lambda_n^l)^{2}-1 }{2}\int
 E_{1}(\lvert u_n^l\rvert^2)\lvert u_n^l\rvert^2dx\notag\\
&\geq
\frac{E(\widetilde{u}_n^l)}{(\lambda_n^l)^{2}}+o(1)\;\;as\;\;
n\rightarrow+\infty,\;\;l\rightarrow+\infty.
\end{align}
From (3.23), we have
\begin{align}
E(\widetilde{U}_n^j)\geq d_m\;\;and\;\;
E(\widetilde{u}_n^l)\geq d_m.
\end{align}
By (3.12), (3.16), (3.25) and (3.26), we deduce that
 as $n\rightarrow +\infty$ and $l\rightarrow+\infty$,
\begin{align}
d_m\geq E(u_n)&=\sum\limits_{j=1}^{l} (\frac{E(\widetilde{U}_n^j)}{\lambda_j^2}+\frac{2(\lambda^{p-1}_j-1)}{p+1}\|U_n^j \|_{L^{p+1}}^{p+1}\notag\\
&+\frac{\lambda_j^{2}-1 }{2}\int E_{1}(\lvert U_n^j\rvert^2)\lvert U_n^j\rvert^2dx)
+\frac{E(\widetilde{u}_n^l)}{(\lambda_n^l)^{2}}+o(1).
\end{align}
Since $1 < p < 3$, combining with (3.15), (3.24) and (3.27), we deduce that by (3.28),
\begin{align}
d_m\geq E(u_n)
%& \geq  \sum\limits_{j=1}^{l}\frac{d_m}{\lambda_j^2}+\inf\limits_{j\geq1}\frac{2(\lambda^{p-1}_j-1)}{p+1}(\sum\limits_{j=1}^{l}\|U_n^j \|_{L^{p+1}(\mathbb{R}^{2})}^{p+1})\notag\\
%&+\inf\limits_{j\geq1}\frac{\lambda_j^{2}-1}{2}(\sum\limits_{j=1}^{l}\int
 %E_{1}(\lvert U_n^j \rvert^2)\lvert U_n^j\rvert^2dx)+\frac{d_m}{(\lambda_n^l)^{2}} +o(1)\notag\\
&\geq  \sum\limits_{j=1}^{l}\frac{d_m}{\lambda_j^2} +\frac{d_m}{(\lambda_n^l)^{2}} + \inf\limits_{j\geq1} (\lambda_j^{p-1}-1)( \frac{1}{2}\int
 E_{1}(\lvert u_n\rvert^2)\lvert u_n\rvert^2dx\notag\\
&+\frac {2}{p+1} \int\lvert u_n\rvert^{p+1}dx)+o(1)\notag\\
&\geq d_m+
 ((\frac{\sqrt{m}}{\|U^{j_0}\|_{L^{2}}})^{p-1}-1)C_0+o(1),
\end{align}
where $C_0>0$ is given in  (3.15). Let $n\rightarrow +\infty$ and $l\rightarrow+\infty$ in (3.29), the following inequality holds
\begin{align}
d_m\geq d_m
+C_0((\frac{\sqrt{m}}{\|U^{j_0}\|_{L^{2}}})^{p-1}-1).
\end{align}
Hence, we get $\|U^{j_0}\|_{L^{2}}^2\geq m$. But by (3.18), we have $\|U^{j_0}\|_{L^{2}}^2\leq m$. Thus $\|U^{j_0}\|_{L^{2}}^2= m$. Put more precisely, in (3.16), there exists only one non-zero term $U^{j_0}$, and the others are zero. Moreover,  from (3.19)-(3.21), it can be obtained that
$E(U^{j_0})=d_m$, and then  the  variational problem (3.7) attains its  infimum at  $U^{j_0}$. Put $Q_{m}=\rvert U^{j_0}\lvert$, which is a minimizer of (3.7).
\par  Step 4. $Q_{m}$ is the positive solution of (2.10).
\par In terms of (3.7), there exists a unique Lagrange multiplier $\omega_{m}$ such that $Q_{m}$ has to satisfy the Euler-Lagrange equation
\begin{align*}
\frac{d}{d\varepsilon}\lvert_{\varepsilon=0}[E(Q_{m}+\varepsilon \eta)+\omega_{m}\int\lvert Q_{m}+\varepsilon \eta\rvert^{2}dx-m\omega_{m}]=0,\;\;\;\;\eta\in C^{\infty}_{0}(\mathbb{R}^{2}).
\end{align*}
It follows that $Q_{m}$ satisfies (2.10) with $\omega=\omega_{m}$. Since $Q_{m}(x)=\rvert U^{j_0}\lvert\geq 0$ a.e in $\mathbb{R}^{2}$, by the strong maximum principle, we get that $Q_{m}(x)=\rvert U^{j_0}\lvert > 0$ for $x\in \mathbb{R}^{2}$. Thus $\rvert U^{j_0}\lvert=Q_{m}(x)$ is a positive minimizer of (3.7). Moreover $Q_{m}$ is the positive solution of (2.10) with $\omega=\omega_{m}$.
\par  This completes the proof of Theorem 3.3.
\par \noindent\textbf{Remark 3.4.} In fact, Ohta \cite{O1995A} solved the variational problem (3.7) with small mass by the concentration compactness principle \cite{L1984A,L1984B}. But here we solve the variational problem (3.7) with definite mass  $0< m < 2d_{J}$ by the profile decomposition.
\par \noindent\textbf{Theorem 3.5.} Let $d_{J}$ be defined as (2.5) and $Q_{\omega}$ be the positive solution of (2.10).
Define
\begin{align}
\mu_{J}=\{\omega\in \mathbb{R}\rvert 0<\int Q^{2}_{\omega}dx < 2d_{J}\}.
\end{align}
Then $\mu_{J}$ is not empty. Moreover $0 < \omega_{J}=\sup\mu_{J} \leq 2d_{J}$.
\par \noindent\textbf{Proof.}  By Theorem 3.3, the Lagrange multiplier $\omega_{m} \in \mu_{J}$. It follows that $\mu_{J}$ is not empty. Then Theorem 2.3 deduces that $0 < \omega_{J}=\sup\mu_{J} \leq 2d_{J}$.
\par  This proves Theorem 3.5.

\par \noindent\textbf{Theorem 3.6.} Let $1 < p < 3$ and $Q_{m}\in H^{1}$ be  a positive minimizer of (3.7). Suppose that the positive solution of (2.10) is unique for every $\omega > 0$. Then the set of all solutions of (3.7) is $S_{m}=\{e^{i\theta}Q_{m}(\cdot+y),\;\;\theta\in \mathbb{R},\;\;y\in \mathbb{R}^{2}\}$. In addition, for arbitrary $u\in S_{m}$, there exists a unique $\omega_{m}>0$ such that $\varphi(t,\;x)=e^{i\omega_{m} t}u(x)$ is a soliton of (1.1).
\par \noindent\textbf{Proof.} From Theorem 3.3, (3.7) has a positive minimizer $Q_{m} \in H^{1}$. Now suppose that $v\in H^{1}$ is an arbitrary solution of (3.7). Let $v=v^{1}+iv^{2}$, where $v^{1},\;v^{2}\in H^{1}$ are real-valued. Then $\widetilde{v}=\lvert v^{1}\rvert+i\lvert v^{2}\rvert$ is still a solution of (3.7). Thus there exists a unique $\omega_{m}>0$ such that $v$ and $\widetilde{v}$ satisfy (2.10). It follows that for $j=1,\;2$,
\begin{align}
\Delta v^{j}+\lvert v\rvert^{p-1}v^{j}+E_{1}(\lvert v\rvert^{2})v^{j}=\omega_{m} v^{j}\;\;in\;\;\mathbb{R}^{2},
\end{align}
\begin{align}
\Delta \lvert v^{j}\rvert+\lvert v\rvert^{p-1}\lvert v^{j}\rvert+E_{1}(\lvert v\rvert^{2})\lvert v^{j}\rvert=\omega_{m} \lvert v^{j}\rvert\;\;in\;\;\mathbb{R}^{2}.
\end{align}
This shows that $\omega$ is the first eigenvalue of the operator $\Delta+\lvert v\rvert^{p-1}+E_{1}(\lvert v\rvert^{2})$ acting over $H^{1}$ and thus, $v^{1},\;v^{2},\;\lvert v^{1}\rvert$ and $\lvert v^{2}\rvert$ are all multiples of a positive normalized eigenfunction $v_{0}$ of $\Delta+\lvert v\rvert^{p-1}+E_{1}(\lvert v\rvert^{2})$, i.e.
\begin{align}
\Delta v_{0}+\lvert v\rvert^{p-1}v_{0}+E_{1}(\lvert v\rvert^{2})v_{0}=\omega v_{0} \;\;in\;\;\mathbb{R}^{2}
\end{align}
with
\begin{align}
v_{0}\in C^{2}(\mathbb{R}^{2})\cap H^{1},\;\;\;\;v_{0}>0\;\;in\;\;\mathbb{R}^{2}\;\;and\;\;\int\lvert v_{0}\rvert^{2}dx=m.
\end{align}
It is now obvious to deduce that: $v=e^{i\theta}v_{0}(\cdot+y)$ for some $\theta\in \mathbb{R}$, $y\in \mathbb{R}^{2}$ and that $v_{0}$ is still a positive solution of (3.7). By the supposition, $v_{0}$ is the unique positive solution of (2.10) with $\omega=\omega_{m}$. It follows that $v_{0}=Q_{m}(\cdot+y)$ for some $y\in \mathbb{R}^{2}$. Thus $v=e^{i\theta}Q_{m}(\cdot+y)$ for some $\theta\in \mathbb{R}$. It is obvious that for any $\theta\in \mathbb{R}$ and $y\in \mathbb{R}^{2}$, $e^{i\theta}Q_{m}(\cdot+y)$ is also a solution of (3.7). Therefore
\begin{align}
S_{m}=\{e^{i\theta}Q_{m}(\cdot+y),\;\;\theta\in \mathbb{R},\;\;y\in \mathbb{R}^{2}\}
\end{align}
is the set of all solutions of (3.7). Moreover for arbitrary $u\in S_{m}$, there exists a unique $\omega_{m}>0$ such that $u$ is a solution of (2.10) with $\omega=\omega_{m}$, which turns out that $\varphi(t,\;x)=e^{i\omega t}u(x)$ is a soliton of (1.1).
\par  This completes the proof of Theorem 3.6.
\par \noindent\textbf{Lemma 3.7} For $1 < p < \infty$ and $u \in H^{1}\backslash \{0\}$, define the functional
\begin{align}
I(u)=2\int \lvert\nabla u\rvert^{2}dx-\frac{2(p-1)}{p+1}\int \lvert u\rvert^{p+1}dx-\int  E_{1}(\lvert u\rvert^{2})\lvert u\rvert^{2}dx.
\end{align}
For $\lambda > 0$, let $u_{\lambda}=\lambda u(\lambda x)$. Then for $\omega > 0$, we have that
\begin{align}
\frac{d}{d\lambda}[E(u_{\lambda})+\omega \int \lvert u_{\lambda}\rvert^{2}dx]=\frac{1}{\lambda}I(u_{\lambda}).
\end{align}
In addition $E(u_{\lambda})+\omega \int \lvert u_{\lambda}\rvert^{2}dx$ attains the minimum at $\lambda_{0}$ satisfying $I(u_{\lambda_{0}})=0$.
Moreover if $u$ is a solution of (2.10), one has that $I(u)=0$.
\par \noindent\textbf{Proof.}  By a direct calculation, it is shown that (3.38) is true. It follows that $E(u_{\lambda})+\omega \int \lvert u_{\lambda}\rvert^{2}dx$ attains the minimum at $\lambda_{0}$ satisfying $I(u_{\lambda_{0}})=0$.
Moreover if $u$ is a solution of (2.10), from (2.11) it follows that $I(u)=0$.
\par  This completes the proof of Lemma 3.7.
\par \noindent\textbf{Theorem 3.8.}  For $1 < p < 3$ and $\omega \in (0, \omega_{J})$, where $\omega_{J}$ is defined as (1.7) and $d_{J}$ is defined as (1.5), we set the constrained variational problem
\begin{align}
d_{\omega}=inf_{\{u\in H^{1}, 0 < \int \rvert u \lvert^{2}dx < 2d_{J}, I(u)=0\}}(E(u)+\omega\int\lvert u\rvert^{2}dx).
\end{align}
Then (3.39) possesses a positive minimizer $Q_{\omega}\in H^{1}$. Moreover $Q_{\omega}$ is the positive solution of  (2.10).
\par \noindent\textbf{Proof.} In the following we complete this proof by five steps.
\par  Step 1. $\{u\in H^{1},  0 < \int \rvert u \lvert^{2}dx < 2d_{J}, I(u)=0\}$ is not empty.
\par  Take $0 < m <2d_{J}$.  From Theorem 3.3 we have that there exists a positive minimizer $Q_{m}(x) \in H^{1}$ such that $ 0 < \int \rvert Q_{m}(x)\lvert^{2}dx < 2d_{J}$ and $Q_{m}(x)$ satisfying (2.10) with $\omega=\omega_{m}\in (0, \omega_{J})$. By Lemma 3.6 it follows that
$I(Q_{m})=0$. Thus $Q_{m}\in\{u\in H^{1},  0 < \int \rvert u \lvert^{2}dx < 2d_{J}, I(u)=0\}$. Therefore $\{u\in H^{1},  0 < \int \rvert u \lvert^{2}dx < 2d_{J}, I(u)=0\}$ is not empty.
\par  Step 2. $d_{\omega} > -\infty$.
\par Take $u\in H^{1}$ satisfying $0<\int \rvert u \lvert^{2}dx < 2d_{J}$ and $I(u)=0$. For $\omega\in (0,\omega_{J})$, we put
\begin{align}
H(u)=E(u)+\omega\int \rvert u \lvert^{2}dx.
\end{align}
From (3.8), it follows that
\begin{align}
 H(u)=&\int\lvert\nabla u\rvert^{2}dx-\frac12\int
 E_{1}(\lvert u\rvert^2)\lvert u\rvert^2dx-\frac {2}{p+1} \int\lvert u\rvert^{p+1}dx+\omega\int \rvert u \lvert^{2}dx\notag\\
&\geq (1-\frac{\|u\|_{L^{2}}^2}{2d_{J}}-\varepsilon)\|\nabla u\|_{L^{2}}^2-C(\varepsilon,\|u\|_{L^{2}}),
\end{align}
where $0<\varepsilon<1-\frac{\|u\|_{L^{2}}^2}{2d_{J}}$. Since $0 < \int \rvert u\lvert^{2}dx < 2d_{J}$, by (3.41)
we have that
\begin{align}
H(u) \geq -C(\varepsilon,\| u\|_{L^{2}})\geq -C(\varepsilon,2d_{J})=constant>-\infty.
\end{align}
\par  Therefore we deduce that $d_{\omega} > -\infty$.
\par  Step 3. Minimizing sequence is bounded in $H^{1}$.
\par  Let $\{u_n\}^{+\infty}_{n=1}\subset \{u\in H^{1}, 0<\int \rvert u \lvert^{2}dx < 2d_{J}, I(u)=0\}$ be a minimizing sequence of (3.39). Then for all $
n\in \mathbb{N}$,
\begin{align}
0 < \int \rvert u_{n} \lvert^{2}dx < 2d_{J},
\end{align}
\begin{align}
H(u_{n})\rightarrow d_{\omega},  n\;\rightarrow\infty.
\end{align}
By (3.41) and (3.44), for $0 < \varepsilon < 1-\frac{1}{2d_{J}}$ and $n$ large enough we have that
\begin{align}
\int \rvert u_{n}\lvert^{2}dx < d_{\omega}+1+C(\varepsilon,2d_{J}).
\end{align}
Combining with (3.41) and (3.43), we deduce that $\{u_n\}^{+\infty}_{n=1}$ is bounded in $H^{1}$.
\par  Step 4. Existence of minimizer.
\par  We apply Lemma 3.1 to the minimizing sequence $\{u_n\}^{+\infty}_{n=1}$. Then there exists a subsequence still denoted by $\{u_n\}^{+\infty}_{n=1}$ such that
\begin{align}
u_n(x)=\sum\limits_{j=1}^lU_n^j(x)+u_n^l,
\end{align}
where $U_n^j(x):=U^j(x-x_n^j)$ and $u_n^l:=u_n^l(x)$ satisfies (3.17). Moreover, as $n\rightarrow +\infty$, (3.18)-(3.21)
are also held. Thus we have
\begin{align}
H(u_n)=\sum\limits_{j=1}^{l}H(U_n^j )+H(u_n^l)+o(1),\;\;as\;n\rightarrow+\infty.
\end{align}
\par  Firstly, we consider the case $d_{\omega} < 0$.
\par  Since $d_{\omega}<0$, by (3.40) and (3.39) for $n$ large enough, we can choose a $0 < \delta < -d_{\omega}$ such that
\begin{align*}
\int\frac12
 E_{1}(\lvert u_n\rvert^2)\lvert u_n\rvert^2+\frac {2}{p+1}\lvert u_n\rvert^{p+1}dx=&\int\lvert \nabla u_n\rvert^{2}+\omega\lvert u_n\rvert^{2}dx-H(u_n)\notag\\
 \geq& -d_\omega-\delta,
\end{align*}
which implies  that for the minimizing sequence $\{u_n\}^{+\infty}_{n=1}$, there exists a constant $C_0>0$ such that for sufficiently large $n$,
\begin{align}
\frac12\int
 E_{1}(\lvert u_n\rvert^2)\lvert u_n\rvert^2dx+\frac {2}{p+1} \int\lvert u_n\rvert^{p+1}dx\geq C_0.
\end{align}
By (3.46), we put $\|u_{n}\|_{L^{2}}^{2}=m$. Then $0<m<2d_{J}$. For $j=1,\cdot\cdot\cdot,l$, let $\widetilde{U}_n^j=\lambda_jU_n^j$ and $\widetilde{u}_n^l=\lambda_n^lu_n^l$,where
\begin{align}
\lambda_j=\frac{\sqrt{m}}{\|U_n^j \|_{L^{2}}}\geq 1,\;\;\;
\lambda_n^l=\frac{\sqrt{m}}{\|u_n^l\|_{L^{2}}}\geq 1.
\end{align}
From the convergence of $\sum\limits_{j=1}^{l}\|U_n^j \|_{L^{2}}^2$, one has that
there exists $j_0\geq 1$ such that
\begin{align}
\inf\limits_{j\geq
1} \lambda_j^{p-1}-1 = \ \lambda_{j_0}^{p-1}-1 = (\frac{\sqrt{m}}{\|U^{j_0}\|_{L^{2}}})^{p-1}-1 .
\end{align}
Now we consider the new energy $H(U_n^j)$ and $H(u_n^l)$. Then we have
\begin{align}
H(U_n^j)=\frac{H(\widetilde{U}_n^j)}{\lambda_j^2}+\frac{2(\lambda_j^{p-1}-1)}{p+1}\int \rvert U_n^j \lvert^{p+1}dx+\frac{  \lambda_j^{2}-1 }{2}\int E_{1}(\lvert U_n^j\rvert^2)\lvert U_n^j\rvert^2dx.
\end{align}
\begin{align}
H(u_n^l)&=\frac{H(\widetilde{u}_n^l)}{(\lambda_n^l)^{2}}+\frac{2(\lambda_n^l)^{p-1}-1)}{p+1}\int \rvert u_n^l\lvert^{p+1}dx+\frac{(\lambda_n^l)^{2}-1 }{2}\int
 E_{1}(\lvert u_n^l\rvert^2)\lvert u_n^l\rvert^2dx\notag\\
&\geq
\frac{H(\widetilde{u}_n^l)}{(\lambda_n^l)^{2}}+o(1),\;\;as\;
n\rightarrow+\infty,\;\;l\rightarrow+\infty.
\end{align}
For $j=1,\cdot\cdot\cdot,l$, let
\begin{align*}
U_{\lambda_j}^j=\mu_{j} \widetilde{U}_{n}^j(\mu_{j} x),\;\;\;v_{n}^l=\mu_{n}^l\widetilde{u}_{n}^l(\mu_{n}^l x).
\end{align*}
Then there exist $0 < \mu_{j}, \mu_{n}^{l}< \infty$ such that
\begin{align*}
I(U_{\lambda_j}^j)=I(v_{n}^l)=0,\;0 < \|U_{\lambda_j}^j \|^{2}_{L^{2}}=\|u_n^l\|^{2}_{L^{2}}=m <  2d_{J}.
\end{align*}
From (3.39), it follows that
\begin{align}
H(U_{\lambda_j}^j)\geq d_{\omega},\;\;H(v_{n}^l)\geq d_{\omega}.
\end{align}
But from Lemma 3.7, one has that
\begin{align}
H(U_{\lambda_j}^j)\leq H(\widetilde{U}_{n}^j),\;\;H(v_{n}^l)\leq H(\widetilde{u}_{n}^l).
\end{align}
Combining with (3.53) and (3.54), we deduce that
\begin{align}
H(\widetilde{U}_{n}^j) \geq d_{\omega},\;\;H(\widetilde{u}_{n}^l) \geq d_{\omega}.
\end{align}
By (3.44), (3.46), (3.51) and (3.52), it can be deduced that
 as $n\rightarrow +\infty$ and $l\rightarrow+\infty$,
\begin{align}
d_m\geq H(u_n)&=\sum\limits_{j=1}^{l} (\frac{H(\widetilde{U}_n^j)}{\lambda_j^2}+\frac{2(\lambda^{p-1}_j-1)}{p+1}\|U_n^j \|_{L^{p+1}(\mathbb{R}^{2})}^{p+1}\notag\\
&+\frac{\lambda_j^{2}-1 }{2}\int E_{1}(\lvert U_n^j\rvert^2)\lvert U_n^j\rvert^2dx)
+\frac{H(\widetilde{u}_n^l)}{(\lambda_n^l)^{2}}+o(1).
\end{align}
Since $1 < p < 3$, combining with (3.48), (3.50), (3.55) and (3.56), we deduce that
\begin{align}
d_m\geq H(u_n)
%& \geq  \sum\limits_{j=1}^{l}\frac{d_m}{\lambda_j^2}+\inf\limits_{j\geq1}\frac{2(\lambda^{p-1}_j-1)}{p+1}(\sum\limits_{j=1}^{l}\|U_n^j \|_{L^{p+1}(\mathbb{R}^{2})}^{p+1})\notag\\
%&+\inf\limits_{j\geq1}\frac{\lambda_j^{2}-1}{2}(\sum\limits_{j=1}^{l}\int
 %E_{1}(\lvert U_n^j \rvert^2)\lvert U_n^j\rvert^2dx)+\frac{d_m}{(\lambda_n^l)^{2}} +o(1)\notag\\
&\geq  \sum\limits_{j=1}^{l}\frac{d_m}{\lambda_j^2} +\frac{d_m}{(\lambda_n^l)^{2}} + \inf\limits_{j\geq1} (\lambda_j^{p-1}-1)( \frac{1}{2}\int
 E_{1}(\lvert u_n\rvert^2)\lvert u_n\rvert^2dx\notag\\
&+\frac {2}{p+1} \int\lvert u_n\rvert^{p+1}dx)+o(1)\notag\\
&\geq d_m+
 ((\frac{\sqrt{m}}{\|U^{j_0}\|_{L^{2}}})^{p-1}-1)C_0+o(1),
\end{align}
where $C_0>0$ is given in  (3.48). Let $n\rightarrow +\infty$ and $l\rightarrow+\infty$ in (3.57), the following inequality holds
\begin{align}
d_m\geq d_m
+C_0((\frac{\sqrt{m}}{\|U^{j_0}\|_{L^{2}}})^{p-1}-1).
\end{align}
Hence, we get   $\|U^{j_0}\|_{L^{2}}^2\geq m$. But by (3.18), we have $\|U^{j_0}\|_{L^{2}}^2\leq m$. Put more precisely, in   (3.48)  there exists only one non-zero term  $U^{j_0}$, and the others are zero. Moreover,  from (3.19)-(3.21), it follows that
$E(U^{j_0})=d_m$, and then  the  variational problem (3.39) attains its  infimum at  $U^{j_0}$.
\par  Secondly, we consider the case $d_{\omega} \geq 0$.
\par  By the profile decomposition, for $n$ large enough, we have
\begin{align}
\sum^{l}_{j=1}H(U_{n}^j)\leq d_{\omega}.
\end{align}
Since $d_{\omega} \geq 0$, there must be some $U_{n}^j$, denoted by $U^j$ such that $H(U^j)\leq d_{\omega}$.
Let $U^{j_{0}}=\lambda U^j(\lambda x)$. There exists $0 < \lambda < \infty$ such that $I(U^{j_{0}})=0$ and $0 < \int \rvert U^{j_{0}}\lvert ^{2}dx < 2d_{J}$. It follows that $ H(U^{j_{0}})\geq d_{\omega}$. Combining with Lemma 3.7, we have $H(U^{j_{0}})= d_{\omega}$
\par  Therefore, no matter $d_{\omega}<0$ or $d_{\omega}\geq0$, there exists $U^{j_{0}}\neq 0$ such that the  variational problem (3.39) attains its  infimum at  $U^{j_0}$. Then $Q_{\omega}(x)=\rvert U^{j_{0}}\lvert \geq 0$ is a minimizer of (3.39).
\par  Step 5. $Q_{\omega}(x)$ is the positive solution of (2.10).
\par  In terms of (3.39), there exists a unique $\Lambda \in \mathbb{R}$ such that $Q_{\omega}(x)=\rvert U^{j_{0}}\lvert$ satisfies the Euler-Lagrange equation for $\eta \in C^{\infty}_{0}(\mathbb{R}^{2})$
\begin{align}
\frac{d}{d\varepsilon}\rvert_{\epsilon =0}(E(Q_{\omega}+\varepsilon\eta)+\omega\int\lvert Q_{\omega}+\varepsilon\eta\rvert^{2}dx+\Lambda I(Q_{\omega}+\varepsilon\eta))=0.
\end{align}
It follows that
\begin{align}
-\Delta Q_{\omega}-Q_{\omega}^{p}- E_{1}(\lvert Q_{\omega}\rvert^{2}) Q_{\omega}+\omega Q_{\omega}+\Lambda(-2\Delta Q_{\omega}-(p-1)Q_{\omega}^{p}-2 E_{1}(\lvert Q_{\omega}\rvert^{2})Q_{\omega})=0.
\end{align}
From (3.61), we have that
\begin{align}
\int (1+2\Lambda)\lvert\nabla Q_{\omega}\rvert^{2}-(1+(p-1)\Lambda)\rvert Q_{\omega}\lvert^{p+1}-(1+2\Lambda) E_{1}(\lvert Q_{\omega}\rvert^{2})\lvert Q_{\omega}\rvert^{2}+\omega \lvert Q_{\omega}\rvert^{2}dx = 0.
\end{align}
\begin{align}
(1+\Lambda(p-1))\frac{2}{p+1}\lvert Q_{\omega}\rvert^{p+1}+\frac{1}{2}(1+2\Lambda)E_{1}(\lvert Q_{\omega}\rvert^{2})\lvert Q_{\omega}\rvert^{2}-\omega \lvert Q_{\omega}\rvert^{2}dx=0.
\end{align}
By $I(Q_{\omega}) = 0$, (3.62) and (3.63), we have that
\begin{align}
\Lambda\int\frac{(p-3)(p-1)}{p+1}\lvert Q_{\omega}\rvert^{p+1}dx=0.
\end{align}
Since $1 < p < 3$ and $Q_{\omega}\geq 0$,  from (3.64), we have $\Lambda=0$. It follows that $Q_{\omega}$ satisfies (2.10). Since
$Q_{\omega}=\lvert U^{j_{0}} \rvert \geq 0$ a.e. in $\mathbb{R}^{2}$, by the strong maximum principle, we get that $Q_{\omega}(x)>0$ for $x\in \mathbb{R}^{2}$. Thus $Q_{\omega}=\lvert U^{j_{0}} \rvert$ is a positive minimizer of (3.39). Moreover $Q_{\omega}$ is the positive solution of (2.10).

\par  This completes the proof of Theorem 3.8.

\par \noindent\textbf{Theorem 3.9.}  Suppose that the positive solution of (2.10) is unique for every $\omega > 0$. Then the variational problem (3.7) determines a one-to-one mapping between $m\in (0, 2d_{J})$ and $\omega\in (0,\omega_{J})$. In detail, for $\omega\in (0,\omega_{J})$ and $m=\int Q^{2}_{\omega}dx$ with positive solution $Q_{\omega}(x)$ of (2.10) , one has that $\frac{dm}{d\omega}=\frac{d}{d\omega}\int Q^{2}_{\omega}dx \neq 0$.
\par \noindent\textbf{Proof.}  For arbitrary $m\in (0, 2d_{J})$, in terms of Theorem 3.3, the variational problem (3.7) determines a positive $Q_{m}(x)\in H^{1}(\mathbb{R}^{2})$ and a unique $\omega_{m}$ such that (2.10) with $\int Q^{2}_{m}dx=m$. By Lemma 2.3 (also see \cite{F2003} and \cite{K2011}), $\omega_{m}\in (0,\omega_{J})$.
\par  Now suppose that there exists another $\omega'\in (0,\omega_{J})$ such that $\omega'\neq \omega_{m}$ and $\int Q^{2}_{\omega'}dx=m$ for the positive solution $Q_{\omega'}(x)$ of (2.10) with $\omega=\omega'$. By the supposition of uniqueness, $Q_{\omega'}(x) > 0$ is unique for (2.10) with $\omega=\omega'$. In addition, $\omega'\neq \omega_{m}$ leads that $Q_{\omega'}(x)\neq Q_{m}(x)$. From Theorem 3.3, $Q_{\omega'}(x)$ is not a minimizer of (3.7). According to Theorem 3.8, $Q_{\omega'}(x)$ must be the positive minimizer of the variational problem (3.39) with $\omega=\omega'$.
\par  We see that $Q_{m}$ satisfies (2.10) with $\omega=\omega_{m}$ and $Q_{\omega'}$ satisfies (2.10) with $\omega=\omega'$. By Lemma 3.7, it follows that $I(Q_{m})=0=I(Q_{\omega'})$.
\par  Summarizing the above facts, we get that
\begin{align}
\int Q^{2}_{\omega'}dx=\int Q^{2}_{m}dx=m;
\end{align}
\begin{align}
I(Q_{\omega'})=I(Q_{m})=0;
\end{align}
\begin{align}
Q_{m}\;is \;the\;minimizer\;of\;(3.7);
\end{align}
\begin{align}
Q_{\omega'}\;is \;the\;minimizer\;of\;(3.39)\;with\;\omega=\omega'.
\end{align}
Since $Q_{\omega'}$ is not a minimizer of (3.7), by (3.65) and (3.67), Theorem 3.3 derives
\begin{align}
E(Q_{m}) < E(Q_{\omega'}).
\end{align}
By (3.65), (3.66) and (3.68), Theorem 3.8 derives that
\begin{align}
E(Q_{\omega'})+\omega'\int Q^{2}_{\omega'}dx \leq E(Q_{m})+\omega'\int Q^{2}_{m}dx.
\end{align}
From (3.65), it is clear that (3.70) is contradictory with (3.69). Therefore it is necessary that $\omega_{m}=\omega'$. It turns that for $\omega \in (0, \omega_{J})$ and positive solution $Q_{\omega}(x)$ of (2.10), we have that
\begin{align*}
\frac{dm}{d\omega}=\frac{d}{d\omega}\int Q^{2}_{\omega}dx\neq 0.
\end{align*}
\par  This completes the proof of Theorem 3.9.
\par \noindent\textbf{Theorem 3.10. }  Let $\mu$ be the set of all Lagrange multipliers corresponding to the all positive minimizers of (3.7).% for $\int q^{2}dx < m < \infty$.
Then $\mu = (0,\omega_{J})$, where $\omega_{J}$ is defined as (1.7). In addition $\mu_{J}=(0,\omega_{J})$, where $\mu_{J}$ is defined as (3.31).
\par \noindent\textbf{Proof.}  Theorem 3.5 derives that $\mu\subset\mu_{J}\subset(0,\omega_{J})$. Now suppose that $\omega \in (0, \omega_{J})$. By Theorem 3.8, this $\omega$ determines a unique $m\in(0,2d_{J})$, and this $m$ determines a unique $\omega_{m}\in \mu$. Then this $\omega_{m}$ can only be $\omega$, that is $\omega_{m}= \omega$. Thus $(0,\omega_{J})\subset \mu$. Therefore $\mu =(0,\omega_{J})$. It follows that $\mu =(0,\omega_{J})=\mu_{J}$.
\par  This proves Theorem 3.10.
\par \noindent\textbf{Remark 3.11. }  Theorem 3.9 shows that for $ m < 2d_{J} $, the normalized solution problem (2.10) with $\int \rvert u\lvert^{2}dx=m$ possesses a unique  $\omega_{m} \in (0, \omega_{J})$ such that $Q_{m}$ is the unique positive solution of (2.10) with $\omega=\omega_{m}$. It gives a positive answer that for (2.10), the mapping from the prescribed mass $m$ to the Lagrange multiplier, that is the soliton frequency $\omega$ is injective. Moreover the approach introduced here can be used to deal with more nonlinear Schr\"odinger equations.

\section{Orbital Stability of Small Solitons}
\par \noindent\textbf{Theorem 4.1.} The soliton $e^{i\omega t}u(x)$ in Theorem 3.6 holds the orbital stability, i.e. for arbitrary $\varepsilon>0$, there exists $\delta>0$ such that for any $\varphi_{0}\in H^{1}$ and $0< \int\lvert\varphi_{0}\rvert^{2}dx < 2d_{J}$, if
\begin{align}
inf_{\{\theta\in \mathbb{R},\;\;y\in \mathbb{R}^{2}\}}\|\varphi_{0}(\cdot)-e^{i\theta}u(\cdot+y)\|_{H^{1}}<\delta,
\end{align}
then the solution $\varphi(t,\;x)$ of the Cauchy problem (1.1)-(2.1) satisfies
\begin{align}
inf_{\{\theta\in \mathbb{R},\;\;y\in \mathbb{R}^{2}\}}\|\varphi(t,\cdot)-e^{i\theta}u(\cdot+y)\|_{H^{1}}<\varepsilon,\;\;t\in \mathbb{R}.
\end{align}
\par \noindent\textbf{Proof.} By Theorem 3.6, it is clear that for arbitrary $u\in S_{m}$ one has that
\begin{align}
S_{m}=\{e^{i\theta}u(\cdot+y),\;\;\theta\in \mathbb{R},\;\;y\in \mathbb{R}^{2}\}.
\end{align}
In terms of Theorem 2.1, for any $\varphi_{0}\in H^{1}$, the Cauchy problem (1.1)-(2.1) possesses a unique global solution $\varphi(t,\;x)\in C(\mathbb{R,}\;H^{1})$ with mass conservation $M(\varphi)=M(\varphi_{0})$ and energy conservation $E(\varphi)=E(\varphi_{0})$ for all $t\in \mathbb{R}$.
\par  Now arguing by contradiction, if the conclusion of Theorem 4.1 does not hold, then there exist $\varepsilon>0$, a sequence $(\varphi^{n}_{0})_{n\in \mathbb{N}^{+}}$ such that
\begin{align}
inf_{\{\theta\in \mathbb{R},\;\;y\in \mathbb{R}^{2}\}}\|\varphi^{n}_{0}-e^{i\theta}u(\cdot+y)\|_{H^{1}}<\frac{1}{n},
\end{align}
and a sequence $(t_{n})_{n\in \mathbb{N}^{+}}$ such that
\begin{align}
inf_{\{\theta\in \mathbb{R},\;\;y\in \mathbb{R}^{2}\}}\|\varphi_{n}(t_{n},\;\cdot)-e^{i\theta}u(\cdot+y)\|_{H^{1}}\geq\varepsilon,
\end{align}
where $\varphi_{n}$ denotes the solution of the Cauchy problem (1.1)-(2.1) with initial datum $\varphi^{n}_{0}$. From (4.4) we yield that
\begin{align}
\int\lvert \varphi^{n}_{0}\rvert^{2}dx\rightarrow\int\lvert u\rvert^{2}dx=m,
\end{align}
\begin{align}
E(\varphi^{n}_{0})\rightarrow E(u)=d_{m}.
\end{align}
Thus (4.6), (4.7), the conservations of mass and energy derive that $\{\varphi_{n}(t_{n},\;\cdot)\}$ is a minimizing sequence for the problem (3.7). Therefore (4.6) and (4.7) derive that there exists $\theta\in \mathbb{R}$ and $y\in \mathbb{R}^{2}$ such that
\begin{align}
\lim_{n\rightarrow \infty}\|\varphi_{n}(t_{n},\;\cdot)-e^{i\theta}u(\cdot+y)\|_{H^{1}}=0.
\end{align}
This is contradictory with (4.5). Theorem 4.1 is proved.
\par \noindent\textbf{Theorem 4.2.} Let $\omega\in (0,\;\infty)$ and $Q_{\omega}$ be a positive solution of (2.10). Then we have
\begin{align}
 \frac{d}{d\omega}E(Q_{\omega})=-\omega\frac{d}{d\omega}M(Q_{\omega})=-\omega\frac{d}{d\omega}\int Q^{2}_{\omega}dx.
\end{align}
\noindent\textbf{Proof.}  Since $Q_{\omega}$ is a positive solution of (2.10), it follows that
\begin{align}
\Delta Q_{\omega}+Q^{p}_{\omega}+ E_{1}(\lvert Q_{\omega}\rvert^{2})Q_{\omega}-\omega Q_{\omega}=0,\;\;\;Q_{\omega}\in H^{1}.
\end{align}
From (2.2) and (4.10), we have
\begin{align*}
\frac{d}{d\omega}E(Q_{\omega})
={} &\frac{d}{d\omega}(\int\lvert\nabla Q_{\omega}\rvert^{2}-\frac{2}{p+1}\lvert Q_{\omega}\rvert^{p+1}-\frac{1}{2} E_{1}(\lvert Q_{\omega}\rvert^{2})\lvert Q_{\omega}\rvert^{2}dx)\notag\\
={} &\int2\lvert \nabla Q_{\omega}\rvert\frac{d}{d\omega}\lvert \nabla Q_{\omega}\rvert-2Q^{p}_{\omega}\frac{d}{d\omega}Q_{\omega}-2 E_{1}(\lvert Q_{\omega}\rvert^{2})Q_{\omega}\frac{d}{d\omega}Q_{\omega}dx\notag\\
={} &\int-2\Delta Q_{\omega}\frac{d}{d\omega}Q_{\omega}-2Q^{p}_{\omega}\frac{d}{d\omega}Q_{\omega}-2 E_{1}(\lvert Q_{\omega}\rvert^{2})Q_{\omega}\frac{d}{d\omega}Q_{\omega}dx\notag\\
={} &\int-2\omega Q_{\omega}\frac{d}{d\omega}Q_{\omega}dx=-\omega\int\frac{d}{d\omega} Q^{2}_{\omega}dx=-\omega\frac{d}{d\omega}\int Q^{2}_{\omega}dx.
\end{align*}
Noting (2.3), this proves (4.9) and completes the proof.
\par  Let $\omega\in (0, \omega_{J})$ and $Q_{\omega}(x)$ be the unique positive solution of (2.10). We set the scalar
\begin{align}
D(\omega)=E(Q_{\omega})+\omega M(Q_{\omega})
\end{align}
and the linearized operator of (4.10)
\begin{align}
H_{\omega}=-\Delta+\omega-pQ^{p-1}_{\omega}-3 E_{1}(\lvert Q_{\omega}\rvert^{2}).
\end{align}
%\begin{align}
%H_{\omega}=-\Delta+\omega-pQ^{p-1}_{\omega}-\frac{1}{2}\mathcal{B}(|Q_{\omega}|^{2})-\frac{1}{2}\mathcal{B}(|\cdot|)|Q_{\omega}|^{2}-2\mathcal{B}(ReQ_{\omega})ReQ_{\omega}.
%\end{align}
It is clear that
\begin{align}
H_{\omega}=\frac{1}{2}E''(Q_{\omega})+\frac{1}{2}\omega M''(Q_{\omega}).
\end{align}
\par \noindent\textbf{Theorem 4.3.}  Let $\omega\in (0, \omega_{J})$ and $Q_{\omega}(x)$ be the unique positive solution of (2.10). Then the operator $H_{\omega}$ has one negative simple eigenvalue and has its kernel spanned by $iQ_{\omega}$. Moreover the positive spectrum of $H_{\omega}$ is bounded away from zero.
\par \noindent\textbf{Proof.} Since $\omega\in (0, \omega_{J})$, By Lemma 2.3, there exists a positive function $Q_{\omega}(x)$ satisfying (4.10). Now suppose that $\lambda\in \mathbb{R}$ satisfies $H_{\omega}Q_{\omega}=\lambda Q_{\omega}$, that is
\begin{align}
-\Delta Q_{\omega}+\omega Q_{\omega}-pQ^{p}_{\omega}-3 E_{1}(\lvert Q_{\omega}\rvert^{2})Q_{\omega}=\lambda Q_{\omega}.
\end{align}
From (4.10), it follows that
\begin{align}
(1-p)Q^{p-1}_{\omega}-2 E_{1}(\lvert Q_{\omega}\rvert^{2})=\lambda.
\end{align}
By Lemma 2.3 and (4.15), we can uniquely determine $\lambda$ as follows
\begin{align}
\lambda=\lambda_{-}=\int(1-p)Q^{p+1}_{\omega}-2 E_{1}(\lvert Q_{\omega}\rvert^{2})\lvert Q_{\omega}\rvert^{2}dx\big/\int Q^{2}_{\omega}dx.
\end{align}
From (4.16), we have $\lambda < 0$.
Therefore we get that $H_{\omega}$ has one negative simple eigenvalue $\lambda_{-}$.  It follows that $H_{\omega}(iQ_{\omega})=\lambda_{-}(iQ_{\omega})$ and the kernel is spanned by $iQ_{\omega}$.
\par  Now suppose that $\lambda > 0$ and $u\in H^{1}\backslash \{0\}$ satisfying $H_{\omega}u=\lambda u$, that is
\begin{align}
-\Delta u+\omega u-pQ^{p-1}_{\omega}u-3E_{1}(\lvert Q_{\omega}\rvert^{2})u=\lambda u.
\end{align}
By Lemma 2.3,
\begin{align}
-pQ^{p-1}_{\omega}-3 E_{1}(\lvert Q_{\omega}\rvert^{2}) : = g(x)=o(\lvert x\rvert^{-1}).
\end{align}
From Kato \cite{K1959}, $-\Delta+g(x)$ has no positive eigenvalues. Thus (4.17) derives that $\lambda \leq \omega$. By Weyl's theorem on the essential spectrum, the rest of the spectrum of $H_{\omega}$ is bounded away from zero (see \cite{RS1978}).
\par  This proves Theorem 4.3.
\par  By Theorem 4.3, $H_{\omega}$ with $T'(0)=i$ satisfies Assumption 3 in \cite{GSS1987} for $\omega\in \mu_{J}$. With $J=-i$, $X=H^{1}$ and $E$ as (2.2), by Theorem 2.1 and Lemma 2.3, (1.1) satisfies Assumption 1 and 2 in \cite{GSS1987} for $\omega\in \mu_{J}$. Thus we can use Theorem 4.7 in \cite{GSS1987} and get the following lemma.
\par \noindent\textbf{Lemma 4.4.} Let $\omega\in (0, \omega_{J})$ and $Q_{\omega}(x)$ be the unique positive solution of (2.10). If $D''(\omega)=\frac{d^{2}}{d\omega^{2}}D(\omega) < 0$, the soliton $e^{i\omega t}Q_{\omega}(x)$ of (1.1) is unstable.
\par  Then we get the following theorem.
\par \noindent\textbf{Theorem 4.5.} Let $\omega\in (0, \omega_{J})$ and $Q_{\omega}(x)$ be the unique positive solution of (2.10). Then we have that
\begin{align*}
\frac{d}{d\omega}\int Q_{\omega}^{2}dx>0.
\end{align*}
\par \noindent\textbf{Proof.} From Theorem 4.2, (4.11) and (4.13), we have that
\begin{align}
D''(\omega)=\frac{d}{d\omega}\int Q_{\omega}^{2}dx.
\end{align}
Since $\omega\in (0, \omega_{J})$, from Theorem 3.10 it follows that $\omega\in \mu_{J}$. In terms of Theorem 4.1, the soliton $e^{i\omega t}Q_{\omega}(x)$ holds the orbital stability. By Lemma 4.4, we deduce that $D''(\omega) \geq 0$. From (4.19) it follows that $\frac{d}{d\omega}\int Q_{\omega}^{2}dx \geq 0$. Set $m(\omega)=\int Q^{2}_{\omega}dx$. From $\omega\in (0, \omega_{J})$, Theorem 3.9 deduces that $\frac{d}{d\omega}\int Q_{\omega}^{2}dx = \frac{dm}{d\omega} \neq 0$. Therefore we get that $\frac{d}{d\omega}\int Q_{\omega}^{2}dx > 0$.
\par  This proves Theorem 4.5.

\par \noindent\textbf{Proof of Theorem A.}  In fact, we have given the proof of Theorem A in the proof of Theorem 4.6. On the other hand, from Theorem 3.5 in \cite{GSS1987}, Theorem 4.5 also deduces Theorem A.

\section{Construction of Multi-Solitons}
\par  It is clear that (1.1) admits the following symmetries.
\par  Time-space translation invariance: if $\varphi(t,x)$ satisfies (1.1), then for any $t_{0},\;x_{0}\in \mathbb{R}\times \mathbb{R}^{2}$,
\begin{align}
\psi(t,x)=\varphi(t-t_{0},x-x_{0})
\end{align}
also satisfies (1.1).
\par  Phase invariance: if $\varphi(t,x)$ satisfies (1.1), then for any $\gamma_{0}\in \mathbb{R}$,
\begin{align}
\psi(t,x)=\varphi(t,x)e^{i\gamma_{0}}
\end{align}
also satisfies (1.1).
\par  Galilean invariance: if $\varphi(t,x)$ satisfies (1.1), then for any $v_{0}\in \mathbb{R}^{2}$,
\begin{align}
\psi(t,x)=\varphi(t,x-v_{0}t)e^{i(\frac{1}{2}v_{0}x-\frac{1}{4}\lvert v_{0}\rvert^{2}t)}
\end{align}
also satisfies (1.1).
\par Let $1 < p < 3$ and $\omega_{J}$ be defined in (1.7). For $K\geq 2$ and $k=1,\;2,\;\cdot\cdot\cdot,\;K$, we take $\omega^{0}_{k}\in (0,\;\omega_{J}),\;\;\gamma^{0}_{k}\in \mathbb{R},\;\;x^{0}_{k}\in \mathbb{R}^{2}$ and $v_{k}\in \mathbb{R}^{2}$ with $v_{k}\neq v_{k'}$ to $k\neq k'$. By Theorem A,
\begin{align}
e^{i\omega^{0}_{k}t}Q_{\omega^{0}_{k}}(x),\;\;k=1,\;2,\;\cdot\cdot\cdot,\;K
\end{align}
are the stable solitons of (1.1). Then in terms of the above symmetries for $k=1,\;2,\;\cdot\cdot\cdot,\;K$,
\begin{eqnarray}
R_{k}(t,x)=Q_{\omega_{k}^{0}}(x-x^{0}_{k}-v_{k}t)e^{i(\frac{1}{2}v_{k}x-\frac{1}{4}\lvert v_{k}\rvert^{2}t+\omega_{k}^{0}t+\gamma_{k}^{0})},\;(t,x)\in\ \mathbb{R}\times \mathbb{R}^{2}
\end{eqnarray}
are also the solitons of (1.1). It is obvious that
\begin{equation}\label{AH4}
\|R_{k}(t,\;\cdot)\|_{L^{2}}=\|Q_{\omega^{0}_{k}}(\cdot)\|_{L^{2}}<\sqrt{2d_{J}},\;t\in\;\mathbb{R},\;k=1,2,\cdot\cdot\cdot,K.
\end{equation}
Now we suppose that $K\geq 2$, $\omega^{0}_{k}\in (0,\;\omega_{J})$ for $k=1,2,\cdot\cdot\cdot,K$, and
\begin{equation}
\sum^{K}_{k=1}\|Q_{\omega^{0}_{k}}(\cdot)\|_{L^{2}}<\sqrt{2d_{J}}.
\end{equation}
Thus
\begin{equation}\label{AH4}
\|\sum^{K}_{k=1}R_{k}(t,\;\cdot)\|_{L^{2}}\leq\sum^{K}_{k=1}\|R_{k}(t,\;\cdot)\|_{L^{2}}<\sqrt{2d_{J}}.
\end{equation}
Now we set
\begin{align}
R(t)=\sum_{k=1}^{K}R_{k}(t,\;\cdot),\;\;t\in \mathbb{R}.
\end{align}
\par \noindent\textbf{Theorem 5.1.} Let $1 < p < 3$. For $K\geq 2$ and $k= {1,\cdot\cdot\cdot,K}$, taking $\omega^{0}_{k}\in (0,\;\omega_{J})$, $\gamma_{k}^{0}\in \mathbb{R},\;x_{k}^{0}\in \mathbb{R}^{2}, \;v_{k}\in \mathbb{R}^{2}$ with $v_{k}\neq v_{k^{'}} \;\;to \;\;k\neq k'$,
$\sum^{K}_{k=1}\|Q_{\omega^{0}_{k}}(\cdot)\|_{L^{2}}<\sqrt{2d_{J}}$ and
\begin{align}
R_{k}(t,x)=Q_{\omega_{k}^{0}}(x-x^{0}_{k}-v_{k}t)e^{i(\frac{1}{2}v_{k}x-\frac{1}{4}\lvert v_{k}\rvert^{2}t+\omega_{k}^{0}t+\gamma_{k}^{0})}
\end{align}
with $(t,x)\in\ \mathbb{R}\times \mathbb{R}^{2}$, there exists a solution $\varphi(t,\;x)$ of (1.1) such that
\begin{align}
\forall t\geq 0,\;\|\varphi(t)-\sum_{k=1}^{K}R_{k}(t)\|_{H^{1}}\leq Ce^{-\theta_{0}t}
\end{align}
for some $\theta_{0}>0$ and $C>0$.
\par \noindent\textbf{Proof of Theorem B.}  Theorem 5.1 directly implies that Theorem B is true.
\par  Let $T_{n}>0$, $n=1,2,\cdot\cdot\cdot$ and $\lim_{n\rightarrow\infty}T_{n}=+\infty$. For $n=1,2,\cdot\cdot\cdot$, by Theorem 2.1 we can let $\varphi_{n}$ be the unique global solution in $H^{1}$ for the Cauchy problem
\begin{equation}\label{AH1}
\left\{
\begin{split}
&i\partial_{t}\varphi_{n}+\Delta \varphi_{n}+\lvert\varphi_{n}\rvert^{p-1}\varphi_{n}+ E_{1}(\lvert\varphi_{n}\rvert^{2})\varphi_{n}=0,&\qquad  (t,x)\in\ \mathbb{R} \times \mathbb{R}^{2},\\
&\varphi_{n}(T_{n},x)=R(T_{n}).
\end{split}
\right.
\end{equation}

 In the following, according to  Martel, Merle and Tsai's way (see \cite{MM2006} and \cite{MMT2006}), we first state the following claim.
\par \noindent\textbf{Claim 5.2.} (Claim 1 in \cite{MM2006}) Let $(v_{k})$, $k=1,\cdot\cdot\cdot,K$ be $K$ vectors of $\mathbb{R}^{2}$ such that for any $k\neq k',v_{k}\neq v_{k'}$.
Then, there exists an orthonormal basis $(e_{1},\;e_{2})$ of $\mathbb{R}^{2}$ such that for any  $k\neq k',(v_{k},e_{1})\neq (v_{k'},e_{1})$.
\par Without any restriction, we can assume that the direction $e_{1}$ given by Claim 5.2 is $x_{1}$, since (1.1) is invariant by rotation. Therefore, we may assume that for any $k\neq k',v_{k,1}\neq v_{k',1}$. We suppose in fact that
\begin{eqnarray}
v_{1,1}<v_{2,1}<\cdot\cdot\cdot<v_{K,1}.
\end{eqnarray}
Since (5.13) and $\omega^{0}_{k}\in (0,\;\omega_{J})$ with $k=1,\cdot\cdot\cdot,K$, we can set $\theta_{0}>0$ such that
\begin{equation}\label{AH4}
\sqrt{\theta_{0}}=\frac{1}{16}min(v_{2,1}-v_{1,1},\cdot\cdot\cdot,v_{K,1}-v_{K-1,1},\;\;\sqrt{\omega_{1}^{0}},\cdot\cdot\cdot,\sqrt{\omega_{K}^{0}}).
\end{equation}
\par Now we state the following uniform estimates about the sequence $(\varphi_{n})$ in (5.12), which is the key point of the proof of Theorem 5.1.
\par \noindent\textbf{Proposition 5.3.} There exist $T_{0}>0,C_{0}>0,\theta_{0}>0$ such that, for all $n \geq 1$,
\begin{equation}\label{AH4}
\forall t\in [T_{0},T_{n}],\;\;\;\|\varphi_{n}(t)-R(t)\|_{H^{1}}\leq C_{0}e^{-\theta_{0}t}.
\end{equation}
\par  In addition, the sequence $(\varphi_{n})$ has the following global bounded property.
\par \noindent\textbf{Lemma 5.4.} There exists a constant $C > 0$, such that, for any $t\in[T_{0},T_{n}]$ and all $n \geq 1$,
\begin{center}
 $\|\varphi_{n}(t)\|_{H^{1}}\leq C$.
 \end{center}
\par \noindent\textbf{Claim 5.5.} ((25) in \cite{CC2011}) Take $\epsilon_{0}>0$. There exists $K_{0}=K_{0}(\epsilon_{0})>0$ such that for all $n$ large enough, we have
\begin{align}
\int_{\lvert x\rvert>K_{0}}\lvert\varphi_{n}(T_{0},x)\rvert^{2}dx \leq \epsilon_{0}.
\end{align}
\par \noindent\textbf{Lemma 5.6.} There exists $\psi_{0}\in H^{1}$ such that up to a subsequence for  $0\leq s<1$
\begin{align}
\varphi_{n}(T_{0})\rightarrow \psi_{0},\;\;\;in\;H^{s}(\mathbb{R}^{2})\;as\;n\rightarrow +\infty.
\end{align}
\par \noindent\textbf{Proof.} By Lemma 5.4, there exists $\psi_{0}\in H^{1}$ such that up to a subsequence,
\begin{align*}
\varphi_{n}(T_{0})\rightharpoonup \psi_{0}\;\;\;in\;H^{1}\;as\;n\rightarrow +\infty.
\end{align*}
From Lemma 5.5, it follows that
\begin{align*}
\varphi_{n}(T_{0})\rightarrow \psi_{0}\;\;\;in\;L^{2}_{loc}(\mathbb{R}^{2})\;as\;n\rightarrow +\infty,
\end{align*}
we conclude that
\begin{align*}
\varphi_{n}(T_{0})\rightarrow \psi_{0}\;\;\;in\;L^{2}\;as\;n\rightarrow +\infty.
\end{align*}
By interpolation we get (5.17).
\par  This completes the proof of Lemma 5.6.
\par \noindent\textbf{Proof of Theorem 5.1.}  Let $\psi_{0}$ be given by Lemma 5.6.
There exists $0<\sigma<1$ such that $1 < p < 1+\frac{4}{2-2\sigma}$ and
\begin{align}
\lvert(\lvert z_{1}\rvert^{p-1}z_{1}+E_{1}(\lvert z_{1}\rvert^{2})z_{1})&-(\lvert z_{2}\rvert^{p-1}z_{2}+E_{1}(\lvert z_{2}\rvert^{2})z_{2})\rvert\notag\\
&\leq C(1+\lvert z_{1}\rvert+\lvert z_{2}\rvert)\lvert z_{1}-z_{2}\rvert\;\;\;
\end{align}
for all $z_{1},\;z_{2}\in \mathbb{C}$.
This implies that the Cauchy problem of (1.1) with $\varphi(T_{0},\;x)=\psi_{0}$ is well-posedness in $H^{\sigma}(\mathbb{R}^{2})$ (see Theorem 5.1.1 in \cite{ C2003}, also refer to \cite{CW1990}). Then we let $\varphi(t,\;x)\in C([T_{0},\;T],\;H^{\sigma}(\mathbb{R}^{2}))$ be the corresponding maximal solution of (1.1) with $\varphi(T_{0},\;x)=\psi_{0}$.  Combining with Lemma 5.6, we can obtain
\begin{align*}
\varphi_{n}(t)\rightarrow \varphi(t)\;\;in\;\;H^{\sigma}(\mathbb{R}^{2})\;\;as\;\;n\rightarrow +\infty
\end{align*}
for any $t\in [T_{0},\;T)$. By boundedness of $\varphi_{n}(t)\;\;in\;\;H^{1}$, we also have
\begin{align*}
\varphi_{n}(t)\rightharpoonup \varphi(t)\;\;in\;\;H^{1}\;\;as\;\;n\rightarrow +\infty
\end{align*}
for any $t\in [T_{0},\;T)$. By Proposition 5.3, for any $t\in [T_{0},\;T)$, we have
\begin{align}
\|\varphi(t)-R(t)\|_{H^{1}}\leq \liminf_{n\rightarrow \infty}\|\varphi_{n}(t)-R(t)\|_{H^{1}}\leq C_{0}e^{-\theta_{0}t}.
\end{align}
In particular, since $R(t)$ is bounded in $H^{1}$ there exists $C>0$ such that for any $t\in [T_{0},\;T)$ we have
\begin{align}
\|\varphi(t)\|_{H^{1}}\leq C_{0}e^{-\theta_{0}t}+\|-R(t)\|_{H^{1}}\leq C.
\end{align}
Recall that, by the blow up alternative (see \cite{C2003}), either $T=+\infty$ or $T<+\infty$ and
$lim_{t\rightarrow T}\|\varphi(t)\|_{H^{1}}=+\infty$. Therefore (5.20) implies that $T=+\infty$. From (5.19) we infer that for all $t\in [T_{0},\;+\infty)$ we have
\begin{align*}
\|\varphi(t)-R(t)\|_{H^{1}}\leq C_{0}e^{-\theta_{0}t}.
\end{align*}
\par  This completes the proof of Theorem 5.1.
\par The proof of the uniform estimates Proposition 5.3 relies on a bootstrap argument. We first state the following bootstrap result.
\par \noindent\textbf{Proposition 5.7.} There exist $A_{0}>0,\theta_{0}>0,T_{0}>0$ and $N_{0}>0$ such that for all $n\geq N_{0}$ and $t^{*}\in [T_{0},T_{n}]$, if
\begin{align}
\forall t\in [t^{*},T_{n}],\;\;\;\|\varphi_{n}(t)-R(t)\|_{H^{1}}\leq A_{0}e^{-\theta_{0}t},
\end{align}
then
\begin{align}
\forall t \in [t^{*},T_{n}],\;\;\;\|\varphi_{n}(t)-R(t)\|_{H^{1}}\leq \frac{A_{0}}{2}e^{-\theta_{0}t}.
\end{align}
By  Proposition 5.7, we deduce the uniform estimates Proposition 5.3.
\par \noindent\textbf{Proof of Proposition 5.3.}(Proposition 1 in \cite{MM2006})
Let $t^{*}$ be the minimal time such that (5.21) holds:
\begin{align*}
t^{*} = min\{\tau\in[T_{0},T_{n}];\;(5.21)\;holds\;for\;all\;t\in[\tau,T_{n}]\}.
\end{align*}
We prove by contradiction that $t^{*}=T_{0}$. Indeed, assume that $t^{*}>T_{0}$. Then
\begin{align*}
\|\varphi_{n}(t^{*})-R(t^{*})\|_{H^{1}}\leq A_{0}e^{-\theta_{0}t},
\end{align*}
and by Proposition 5.7 we can improve this estimate in
\begin{align*}
\|\varphi_{n}(t^{*})-R(t^{*})\|_{H^{1}}\leq \frac{A_{0}}{2}e^{-\theta_{0}t}.
\end{align*}
Hence, by continuity of $\varphi_{n}(t)$ in $H^{1}$, there exists $T_{0}\leq t^{**}<t^{*}$ such that (5.21) holds for all $t\in [t^{**},\;t^{*}]$. This contradicts the minimality of $t^{*}$.
\par  This completes the proof of Proposition 5.3.
\par  Now for $k=1,\;\cdot\cdot\cdot,\;K$, let $\omega_{k}\in (0,\;\omega_{J})$ and $Q_{\omega_{k}}(x)$ be the unique positive solutions of (2.10). To $ x^{0}_{k},\;x_{k},\;v_{k}\in \mathbb{R}^{d}$ and $\gamma_{k}\in \mathbb{R},\;k=1,\cdot\cdot\cdot,K$, we assume that
\begin{align*}
\widetilde{R}_{k}=Q_{\omega_{k}}(\cdot-\widetilde{x}_{k})e^{i(\frac{1}{2}v_{k}x+\delta_{k})},\;\;\widetilde{x}_{k}=x^{0}_{k}+v_{k}t+x_{k},\;\;\delta_{k}=-\frac{1}{4}\rvert v_{k}\lvert^{2}t+\omega^{0}_{k}t+\gamma_{k},
\end{align*}
\begin{align*}
\widetilde{R}=\sum^{K}_{k=1}\widetilde{R}_{k}\;\;\;\;and\;\;\;\;\;\;\varepsilon=\varphi_{n}-\widetilde{R}.
\end{align*}
\par  For $\alpha>0$, $l>0$, $\omega^{0}_{k} \in (0,\;\omega_{J})$, $\widetilde{\gamma}_{k}\in \mathbb{R}$ and $\widetilde{y}_{k}\in \mathbb{R}^{2}$, $k=1,\cdot\cdot\cdot,K$ set
\begin{align}
&\mu(\alpha,l)=\{\varphi_{n}\in H^{1};\;\notag\\
&inf_{\{\widetilde{\gamma}_{k}\in \mathbb{R},\rvert \widetilde{y}_{k}\lvert-\rvert\widetilde{y}_{k-1}\lvert>l\}}\|\varphi_{n}(t,\;\cdot)-\sum^{K}_{k=1}Q_{\omega^{0}_{k}}(\cdot-\widetilde{y}_{k})e^{i(\frac{1}{2}v_{k}x+\widetilde{\gamma}_{k})}\|_{H^{1}}<\alpha\}.
\end{align}
%and $\alpha=C_{1}A_{0}e^{-\theta_{0}t}$ \\ \\
\par \noindent\textbf{Lemma 5.8.} There exists $\alpha_{1}>0$, $C_{1}>0$, $l_{1}>0$, and a unique $C^{1}$ function $(\omega_{k},x_{k},\gamma_{k}):\mu(\alpha_{1},\;l_{1})\rightarrow (0,\omega_{J})\times \mathbb{R}^{2}\times \mathbb{R}$ for any $k=1,\cdot\cdot\cdot,K$, such that if $\varphi_{n}\in \mu(\alpha_{1},\;l_{1})$, then
\begin{align}
Re\int \widetilde{R}_{k}\overline{\varepsilon}dx=Im\int \widetilde{R}_{k}\overline{\varepsilon}dx=0,\;\;\;Re\int \nabla Q_{\omega_{k}}(\cdot-\widetilde{x}_{k})e^{i(\frac{1}{2}v_{k}x+\delta_{k})}\overline{\varepsilon}dx=0.
\end{align}
Moreover, if $\varphi_{n}\in \mu(\alpha, l)$, for $0<\alpha<\alpha_{1}$, $0<l_{1}<l$, then
\begin{equation}
\|\varepsilon\|_{H^{1}}+\sum^{K}_{k=1}\rvert\omega_{k}-\omega^{0}_{k}\lvert\leq C_{1}\alpha,\;\;\rvert\widetilde{x}_{k}\lvert-\rvert\widetilde{x}_{k-1}\lvert>l-C_{1}\alpha>\frac{l}{2}.
\end{equation}
\par \noindent\textbf{Proof.}  The proof is a standard application of the implicit function. Let $\alpha > 0$ and $L > 0$.
Let $\omega^{0}_{1},\cdot\cdot\cdot,\omega^{0}_{K}\in (0,\omega_{J})$, $\gamma^{0}_{1},\cdot\cdot\cdot,\gamma^{0}_{K}\in \mathbb{R}$, and $\widetilde{x}^{0}_{1},\cdot\cdot\cdot,\widetilde{x}^{0}_{K}\in \mathbb{R}^{2}$ such that $\rvert\widetilde{x}^{0}_{k}\lvert>\rvert\widetilde{x}^{0}_{k-1}\lvert+l$. Let $B_{0}$ be the $B_{0}-ball$ of center $\sum^{K}_{k=1}R_{k}$ with $R_{k}=Q_{\omega^{0}_{k}}(\cdot-\widetilde{x}^{0}_{k})e^{i(\frac{1}{2}v_{k}x+\delta^{0}_{k})}$, where
$\widetilde{x}^{0}_{k}=x^{0}_{k}+v_{k}t,\;\;\delta^{0}_{k}=-\frac{1}{4}\rvert v_{k}\lvert^{2}t+\omega^{0}_{k}t+\gamma^{0}_{k}$
and of radius $10\alpha$. For any $\varphi_{n}\in B_{0}$ and parameters $\omega_{1},\cdot\cdot\cdot,\omega_{K};\widetilde{x}_{1},\cdot\cdot\cdot,\widetilde{x}_{K};\gamma_{1},\cdot\cdot\cdot,\gamma_{K}$, let
$s=(\omega_{1},\cdot\cdot\cdot,\omega_{K};\widetilde{x}_{1},\cdot\cdot\cdot,\widetilde{x}_{K};\gamma_{1},\cdot\cdot\cdot,\gamma_{K};\varphi_{n})$.
Define the following functions of $s$
\begin{eqnarray*}
\rho^{1}_{k}(s)=Re\int \widetilde{R}_{k}\overline{\varepsilon}(s;x)dx;\;\;\rho^{2}_{k}(s)=Re\int\nabla Q_{\omega_{k}}(\cdot-\widetilde{x}_{k})e^{i(\frac{1}{2}v_{k}x+\delta_{k})}\overline{\varepsilon}(s;x)dx;\;\;
\end{eqnarray*}
\begin{align*}
\rho^{3}_{k}(s)=Im\int \widetilde{R}_{k}\overline{\varepsilon}(s;x)dx,
\end{align*}
for $s$ close to $s_{0}=(\omega^{0}_{1},\cdot\cdot\cdot,\omega^{0}_{K};\widetilde{x}^{0}_{1},\cdot\cdot\cdot,\widetilde{x}^{0}_{K};
\gamma^{0}_{1},\cdot\cdot\cdot,\gamma^{0}_{K};\sum^{K}_{k=1}R_{k})$.
\par  When $s=s_{0}$, we have $\varepsilon(s_{0})= 0$, and thus for $j=1,2,3$, $\rho^{j}_{k}(s_{0})=0$. For $\varphi_{n}\in B_{0}$, we can apply the implicit theorem to prove (5.24). It means that we can choose the unique coefficients $(\omega_{1},\cdot\cdot\cdot,\omega_{K};\widetilde{x}_{1},\cdot\cdot\cdot,\widetilde{x}_{K};\gamma_{1},\cdot\cdot\cdot,\gamma_{K})$, such that $s$ is close to $s_{0}$ and verifies $\rho^{j}_{k}(s)=0$ for $j=1,2,3$. In order to apply the implicit function theorem to this situation, we compute the derivatives of $\rho^{j}_{k}$ for any $k,j$ corresponding to each $(\omega_{k},\widetilde{x}_{k},\gamma_{k})$. Note that
\begin{align*}
\frac{\partial\varepsilon}{\partial\omega_{k}}(s_{0})=-\frac{\partial Q_{\omega}}{\partial\omega}\bigg\rvert_{\omega=\omega^{0}_{k}}
(\cdot-\widetilde{x}^{0}_{k})e^{i(\frac{1}{2}v_{k}x+\delta^{0}_{k})},
\end{align*}
\begin{align*}
\nabla_{x_{k}}\varepsilon(s_{0})=\nabla  Q_{\omega^{0}_{k}}(\cdot-\widetilde{x}^{0}_{k})e^{i(\frac{1}{2}v_{k}x+\delta^{0}_{k})},\;\;\;\;\;\frac{\partial\varepsilon}{\partial \gamma_{k}}(s_{0})=-iR_{k}.
%\frac{\partial\varepsilon}{\partialx_{k}}(q_{0})=Q'_{\omega^{0}_{k}}(\cdot-x^{0}_{k}-v_{k}t)e^{i(\frac{1}{2}v_{k}x-\frac{1}{4}|v_{k}|^{2}t+\omega^{0}_{k}t+\gamma^{0}_{k})},
\end{align*}
Thus for $j=1$
\begin{align*}
\frac{\partial\rho^{1}_{k'}}{\partial\omega_{k}}(s_{0})=-Re \int R_{k'}\frac{\partial Q_{\omega}}{\partial\omega}\bigg\rvert_{\omega=\omega^{0}_{k}}
(\cdot-\widetilde{x}^{0}_{k})e^{-i(\frac{1}{2}v_{k}x+\delta^{0}_{k})}dx,
\end{align*}
\begin{align*}
\nabla_{\widetilde{x}_{k}}\rho^{1}_{k'}(s_{0})=Re \int R_{k'}\nabla  Q_{\omega^{0}_{k}}(\cdot-\widetilde{x}^{0}_{k})e^{-i(\frac{1}{2}v_{k}x+\delta^{0}_{k})}dx,\;\;\;\;\frac{\partial\rho^{1}_{k'}}{\partial\gamma_{k}}(s_{0})=-Im \int R_{k'} \overline{R}_{k}dx,
\end{align*}
and similar formulas hold for $\frac{\partial\rho^{2}_{k'}}{\partial\omega_{k}}(s_{0})$, $\frac{\partial\rho^{2}_{k'}}{\partial x_{k}}(s_{0})$, $\frac{\partial\rho^{2}_{k'}}{\partial\gamma_{k}}(s_{0})$, $\frac{\partial\rho^{3}_{k'}}{\partial\omega_{k}}(s_{0})$, $\frac{\partial\rho^{3}_{k'}}{\partial x_{k}}(s_{0})$ and $\frac{\partial\rho^{3}_{k'}}{\partial\gamma_{k}}(s_{0})$. For $k'=k$, by Theorem 4.6, we have
\begin{align}
\frac{\partial\rho^{1}_{k}}{\partial \omega_{k}}(s_{0})=a_{k}<0,\;\;\;
\frac{\partial\rho^{2}_{k}}{\partial \omega_{k}}(s_{0})=0,\;\;\;
\frac{\partial\rho^{3}_{k}}{\partial \omega_{k}}(s_{0})=0;\;\;\;
\end{align}
\begin{align}
\nabla_{\widetilde{x}_{k}}\rho^{1}_{k}(s_{0})=0,\;\;\;
\nabla_{\widetilde{x}_{k}}\rho^{2}_{k}(s_{0})=b_{k}>0,\;\;\;
\nabla_{\widetilde{x}_{k}}\rho^{3}_{k}(s_{0})=0;\;\;\;
\end{align}
\begin{align}
\frac{\partial\rho^{1}_{k}}{\partial \gamma_{k}}(s_{0})=0,\;\;\;
\frac{\partial\rho^{2}_{k}}{\partial \gamma_{k}}(s_{0})=0,\;\;\;
\frac{\partial\rho^{3}_{k}}{\partial \gamma_{k}}(s_{0})=c_{k}>0.\;\;\;
\end{align}
For $k'\neq k$ and $j=1,2,3$, by Lemma 2.4, we know the different $Q_{\omega_{k}}$ are exponentially decaying and located at centers distant at least of $l$, thus we have
\begin{equation}
\big\rvert\frac{\partial\rho^{j}_{k'}}{\partial \omega_{k}}(s_{0})\big\lvert+
\big\rvert\nabla_{\widetilde{x}_{k}}\rho^{j}_{k'}(s_{0})\big\lvert+
\big\rvert\frac{\partial\rho^{j}_{k'}}{\partial \gamma_{k}}(s_{0})\big\lvert\leq Ce^{-\theta_{0}l}.
\end{equation}
These terms are arbitrarily small by choosing $l$ large enough.
\par    By (5.26), (5.27), (5.28) and (5.29), we know the Jacobian of $\rho=(\rho^{1}_{1},\cdot\cdot\cdot,\rho^{1}_{K};\rho^{2}_{1},\cdot\cdot\cdot,\rho^{2}_{K};\rho^{3}_{1},\cdot\cdot\cdot,\rho^{3}_{K})$
as a function of
$(\omega_{1},\cdot\cdot\cdot,\omega_{K};\widetilde{x}_{1},\cdot\cdot\cdot,\widetilde{x}_{K};\gamma_{1},\cdot\cdot\cdot,\gamma_{K})$ at the point $s_{0}$ is not zero. By the implicit function theorem, for $\alpha$ small and $\varphi_{n}\in B_{0}$, there exist unique parameters $(\omega_{1},\cdot\cdot\cdot,\omega_{K};\widetilde{x}_{1},\cdot\cdot\cdot,\widetilde{x}_{K};\gamma_{1},\cdot\cdot\cdot,\gamma_{K})$ such that $\rho(s)=0$. We obtain directly estimates (5.24) with constants that are independent of the ball $B_{0}$. This proves the result for $\varphi_{n}\in B_{0}$. If we now take $\varphi_{n}\in \mu (\alpha, l)$, then $\varphi_{n}\in $ belongs to such a ball $B_{0}$, and the results follows.
\par  This completes the proof of Lemma 5.8.
\par  By Lemma 5.8, we see that $\omega_{k}$, $\gamma_{k}$ and $x_{k}$  are all functions of $t\in [t^{*},T_{n}]$, that is $\omega_{k}=\omega_{k}(t)$, $\gamma_{k}=\gamma_{k}(t)$ and $x_{k}=x_{k}(t)$. Thus we replace the former assumptions about $\widetilde{R}_{k}$, $\widetilde{R}$ and $\varepsilon$ as follows.
\par  For $k=1,\;\cdot\cdot\cdot,\;K$, let $\omega_{k}(t)\in (0,\;\omega_{J})$ and $Q_{\omega_{k}(t)}(x)$ be the positive solutions of (2.10). To $ x^{0}_{k},\;x_{k}(t),\;v_{k}\in \mathbb{R}^{2}$ and $\gamma_{k}(t)\in \mathbb{R},\;k=1,\cdot\cdot\cdot,K$, we set $\widetilde{x}_{k}(t)=x^{0}_{k}+v_{k}t+x_{k}(t),\;\;\delta_{k}(t)=-\frac{1}{4}\lvert v_{k}\rvert^{2}t+\omega^{0}_{k}t+\gamma_{k}(t)$,
\begin{align}
\widetilde{R}_{k}(t)=Q_{\omega_{k}(t)}(\cdot-\widetilde{x}_{k}(t))e^{i(\frac{1}{2}v_{k}x+\delta_{k}(t))},\;\;
\end{align}
\begin{eqnarray}
\widetilde{R}(t)=\sum^{K}_{k=1}\widetilde{R}_{k}(t)\;\;\;\;and\;\;\;\;\;\;\varepsilon(t,\cdot)=\varphi_{n}(t,\cdot)-\widetilde{R}(t).
\end{eqnarray}
\par \noindent\textbf{Lemma 5.9.} (Lemma 3 in \cite{MM2006}) There exists $C_{1}>0$ such that if $T_{0}$ is large enough, then there exists a unique $C^{1}$ function $(\omega_{k},x_{k},\gamma_{k}):[t^{*},T_{n}]\rightarrow(0,\;\omega_{J})\times \mathbb{R}^{2}\times \mathbb{R}$, for any $k=1,2,\cdot\cdot\cdot,K$ such that
\begin{equation}\label{AH4}
Re\int \widetilde{R}_{k}(t)\overline{\varepsilon}(t)dx=Im\int \widetilde{R}_{k}(t)\overline{\varepsilon}(t)dx=0,\;\;Re\int \nabla \widetilde{R}_{k}(t)\overline{\varepsilon}(t)dx=0,
\end{equation}
\begin{equation}\label{AH4}
\|\varepsilon(t)\|_{H^{1}}+\sum_{k=1}^{K}\lvert\omega_{k}(t)-\omega_{k}^{0}\rvert\leq C_{1}A_{0}e^{-\theta_{0}t},
\end{equation}
and
\begin{equation}\label{AH4}
\lvert\dot{\omega}_{k}(t)\rvert^{2}+\lvert\dot{x}_{k}(t)\rvert^{2}+\lvert\dot{\gamma}_{k}(t)-(\omega_{k}(t)-\omega_{k}^{0})\rvert^{2}\leq C_{1}\|\varepsilon(t)\|^{2}_{H^{1}}+C_{1}e^{-2\theta_{0}t}.
\end{equation}
\noindent\textbf{Proof.}  The first part of the statement follows from Lemma 5.8, hence the main thing to check is (5.34). We first write the equation verified by $\varepsilon$. Recall that $\varphi_{n}$ satisfies $i\partial_{t}\varphi_{n}=E'(\varphi_{n})$, we replace $\varphi_{n}$ by $\varepsilon(t)+\widetilde{R}(t)$ in the previous equation to get
\begin{align}
i\partial_{t}\varepsilon+\mathcal{L}(\varepsilon)=
&-i\sum_{k=1}^{K}[\dot{\omega}_{k}(t)\frac{\partial Q_{\omega}}{\partial\omega}\bigg\rvert_{\omega=\omega_{k}(t)}(\cdot-\widetilde{x}_{k}(t))e^{i(\frac{1}{2}v_{k}x+\delta_{k}(t))}]\notag\\
&+i\sum_{k=1}^{K}[\dot{x}_{k}(t)\nabla Q_{\omega_{k}(t)}(\cdot-\widetilde{x}_{k}(t))e^{i(\frac{1}{2}v_{k}x+\delta_{k}(t))}]\notag\\
&+\sum_{k=1}^{K}[(\dot{\gamma_{k}}(t)-(\omega_{k}(t)-\omega^{0}_{k}))\widetilde{R}_{k}(t)]+\mathcal{N}(\varepsilon)+O(e^{-2\theta_{0}t}),
\end{align}
where
\begin{align*}
\mathcal{L}(\varepsilon):=&\Delta \varepsilon+\sum^{K}_{k=1}(\rvert \widetilde{R}_{k}\lvert^{p-1}\varepsilon +E_{1}(\rvert \widetilde{R}_{k}\lvert^{2})\varepsilon\notag\\
+&((p-1)\rvert \widetilde{R}_{k}\lvert^{p-3}+2E_{1}(\rvert \widetilde{R}_{k}\lvert^{2}))Re(\widetilde{R}_{k}\overline{\varepsilon})\widetilde{R}_{k}),
\end{align*}
and $\mathcal{N}(\varepsilon)$ is the remaining nonlinear part.
\par  Now take the scalar product of (5.35) with $i\widetilde{R}_{k}$, $\widetilde{R}_{k}$, $\partial_{x}\widetilde{R}_{k}$. By the definition of $\widetilde{R}_{k}$, exponential localization and the orthogonality condition (5.32), we obtain a differential system for the modulation equations vector $Mod(t)=(\dot{\omega}_{k}(t), \dot{x}_{k}(t), \dot{\gamma_{k}}(t)-(\omega_{k}(t)-\omega^{0}_{k})),\;k=1,2,\cdot\cdot\cdot,K$ of the form
\begin{align}
Mod(t)=B(\varepsilon)+O(e^{-2\theta_{0}t}),
\end{align}
where $\rvert B(\varepsilon) \lvert \leq M\|\varepsilon\|_{H^{1}}$. As long as the modulation parameter do not vary too much and $\|\varepsilon\|_{H^{1}}$ remains small, $M$ is invertible and we can deduce that
\begin{align}
\rvert Mod(t)\lvert\leq M\|\varepsilon\|_{H^{1}}+O(e^{-2\theta_{0}t}).
\end{align}
Thus one deduces that (5.34).
\par This completes the proof of Lemma 5.9.
\par \noindent\textbf{Claim 5.10.}(Claim 2 in \cite{MM2006}) Let $z(t)\in H^{1}$ be a solution of (1.1). Let $h:x_{1}\in \mathbb{R}\mapsto h(x_{1})$ be a $C^{3}$ real-valued function of one variable such that $h$, $h'$ and $h'''$ are bounded. Then, for all $t\in \mathbb{R}$
\begin{eqnarray}
\frac{1}{2}\frac{d}{dt}\int\lvert z\rvert^{2}h(x_{1})dx=Im\int\partial_{x_{1}}z\bar zh'(x_{1})dx,
\end{eqnarray}
\begin{align}
\frac{1}{2}\frac{d}{dt}Im\int\partial_{x_{1}}z\bar zh(x_{1})dx=&\int\lvert \partial_{x_{1}}z\rvert^{2}h'(x_{1})dx-\frac{p-1}{2(p+1)}\int\lvert z\rvert^{p+1}h'(x_{1})dx\notag\\
-&\frac{1}{4}\int\lvert z\rvert^{2}h'''(x_{1})dx+\frac{1}{4}\int\lvert \nabla z_{n}\rvert^{2}h'(x_{1})dx\notag\\
-&\frac{1}{2}\int E_{1}(\lvert z\rvert^{2})\lvert z\rvert^{2}h'(x_{1})dx,
\end{align}
where $\partial_{x_{1}}z_{n}= E_{1}(\lvert z\rvert^{2})$ and $\Delta z_{n}=\partial_{x_{1}}\rvert z\lvert^{2}$.
\begin{align}
\frac{1}{2}\frac{d}{dt}Im\int\partial_{x_{2}}z\bar zh(x_{1})dx&=Re\int\partial_{x_{2}}z\partial_{x_{1}}\bar zh'(x_{1})dx
+\frac{1}{2}\int\partial_{x_{1}}z_{n}\partial_{x_{2}}z_{n}h'(x_{1})dx\notag\\
&-\frac{1}{2}\int\partial_{x_{2}}z_{n}\cdot \rvert z\lvert^{2}h'(x_{1})dx.
\end{align}
Since $\varphi_{n}(T_{n})=R(T_{n})$ and at time $t=T_{n}$  the decomposition in (5.24) is unique, it follows that
\begin{equation}\label{AH4}
\varepsilon(T_{n})\equiv 0,\;\;\widetilde{R}(T_{n})\equiv R(T_{n}),\;\;\omega_{k}(T_{n})=\omega_{k}^{0},\;\;x_{k}(T_{n})=0,\;\;\gamma_{k}(T_{n})=\gamma_{k}^{0}.
\end{equation}
\par  Let $Y(s)$ be a $C^{3}$ function such that
\begin{equation}\label{AH4}
0\leq Y\leq 1\;\;on\;\mathbb{R};\;\;\;Y(s)=0\;\;for\;s\leq -1;\;\;\;Y(s)=1\;\;for\;s>1;\;\;\;Y'\geq 0\;\;\;on\;\mathbb{R}
\end{equation}
and satisfying for some constant $C>0$,
\begin{eqnarray*}
(Y'(x))^{2}\leq C Y(x),\;\;\;(Y''(x))^{2}\leq C Y'(x)\;\;\;for\;all\;x\in \mathbb{R}.
\end{eqnarray*}
For this, consider $Y(s)=\frac{1}{16}(1+s)^{4}\;\;for\;\;s \in (-1,\;0)$ close to $-1$, and similarly at $s=1$.
\par  For all $k=2,\cdot\cdot\cdot,K$, let
\begin{eqnarray*}
\sigma_{k}=\frac{1}{2}(v_{k-1,1}+v_{k,1}).
\end{eqnarray*}
For $L>0$ large enough to be fixed later, for any $k=2,\cdot\cdot\cdot,K-1$, let
\begin{eqnarray}
y_{k}(t,x)=Y(\frac{x_{1}-\sigma_{k}{t}}{L})-Y(\frac{x_{1}-\sigma_{k+1}{t}}{L}),
\end{eqnarray}
\begin{equation}
y_{1}(t,x)=1-Y(\frac{x_{1}-\sigma_{2}{t}}{L}),\;\;y_{K}(t,x)=Y(\frac{x_{1}-\sigma_{K}{t}}{L}).
\end{equation}
Finally, set for all $k=1,\;\cdot\cdot\cdot,\;K$:
\begin{equation}
I_{k}(t)=\int\lvert\varphi_{n}(t,x)\rvert^{2}y_{k}(t,x)dx,\;\;\;M_{k}(t)=Im \int \nabla \varphi_{n}(t,x)\bar \varphi_{n}(t,x)y_{k}(t,x)dx.
\end{equation}
The quantities $I_{k}(t)$ and $M_{k}(t)$ are local versions of the $L^{2}$ norm and momentum. Ordering the $v_{k,\;1}$ as in (5.13) was useful to split the various solitons using only the coordinate $x_{1}$.
\par \noindent\textbf{Lemma 5.11.}(Lemma 3.5 in \cite{WC2017}) Let $L>0$. There exists $C>0$ such that if $L$ and $T_{0}$ are large enough, then for all $k=2,\;\cdot\cdot\cdot,\;K$, $t\in[t^{*},T_{n}]$, we have
\begin{eqnarray}
\lvert I_{k}(T_{n})-I_{k}(t)\rvert+\lvert M_{k}(T_{n})-M_{k}(t)\rvert\leq \frac{CA_{0}^{2}}{L}e^{-2\theta_{0}t}.%\;\;\;k=2,\cdot\cdot\cdot,K.
\end{eqnarray}
\noindent\textbf{Proof.} From (5.38), we have
\begin{align}
\frac{1}{2}\frac{d}{dt}\int\rvert\varphi_{n}\lvert^{2}Ydx=\frac{1}{L}Im\int\partial_{x_{1}}\varphi_{n}\bar \varphi_{n} Y'dx
-\frac{\sigma_{k}}{2L}\int\rvert\varphi_{n}\lvert^{2}Y'dx.
\end{align}
Set
\begin{eqnarray*}
\Omega_{1}=\Omega_{1}(t)=[-L+\sigma_{k}t,L+\sigma_{k}t]\times \mathbb{R}.
\end{eqnarray*}
Thus, by the properties of $Y$ and (5.47), we obtain
\begin{eqnarray}
\rvert\frac{d}{dt}\int\rvert\varphi_{n}\lvert^{2}Ydx\lvert\leq \frac{C}{L}\int_{\Omega_{1}}(\rvert\partial_{x_{1}}\varphi_{n}\lvert^{2}+\rvert\varphi_{n}\lvert^{2})dx. %\frac{C}{L}\int_{\Omega_{1}}(|\nabla\varphi_{n}|^{2}+|\varphi_{n}|^{2})dx.
\end{eqnarray}
Similarly, by (5.39), we have
\begin{align}
\frac{1}{2}\frac{d}{dt}Im\int \partial_{x_{1}}\varphi_{n}\bar \varphi_{n} Ydx
{}&=\frac{1}{L}\int(\rvert\partial_{x_{1}}\varphi_{n}\lvert^{2}
-\frac{p-1}{2(p+1)}\rvert\varphi_{n}\lvert^{p+1})Y'dx\notag\\
&-\frac{1}{4L^{3}}\int\rvert\varphi_{n}\lvert^{2}Y'''dx-\frac{\sigma_{k}}{2L}Im\int\partial_{x_{1}}\varphi_{n}\bar \varphi_{n} Y'dx\notag\\
&+\frac{1}{L}\int (-E_{1}(\rvert\varphi_{n}\lvert^{2})\rvert\varphi_{n}\lvert^{2}+\frac{1}{2}\rvert \nabla z_{n}\lvert^{2})Y'dx.
\end{align}
Notice that $\nabla z_{n}=(E_{1}(\rvert \varphi_{n}\lvert^{2}), E_{2}(\rvert \varphi_{n}\lvert^{2}))$. To obtain time decay of the variation of momentum, we decompose $\varphi_{n}=\sum^{K}_{k=1}R_{k}+\varepsilon$ to obtain
\begin{align}
\int& _{\Omega_{1}}\rvert E_{1}(\rvert \varphi_{n}\lvert^{2})\lvert \rvert \varphi_{n}\lvert^{2}dx=\notag\\
&\int _{\Omega_{1}}\{\sum^{K}_{k=1}\rvert E_{1}(\rvert R_{k}\lvert^{2})\lvert+2\sum_{k\neq k'}\rvert E_{1}(Re(R_{k}\overline{R}_{k})) \lvert+\rvert E_{1}(\rvert \varepsilon\lvert^{2})\lvert\}\{\sum^{K}_{k=1}\rvert R_{k}\lvert^{2}+2Re(R_{k}\overline{R}_{k})+\rvert \varepsilon\lvert^{2}\}dx
\end{align}
and
\begin{align}
\int _{\Omega_{1}}\rvert z_{n}\lvert^{2}dx=\sum^{2}_{n=1}\int _{\Omega_{1}}\{\sum^{K}_{k=1}\rvert E_{n}(\rvert R_{k}\lvert^{2})\lvert+2\sum_{k\neq k'}\rvert E_{n}(Re(R_{k}\overline{R}_{k'})) \lvert+\rvert E_{n}(\rvert \varepsilon\lvert^{2})\lvert\}^{2}dx.
\end{align}
By Lemma 2.4, we estimate each term of (5.50) and (5.51) separately as follows:
\begin{align*}
\int _{\Omega_{1}}E_{1}(\rvert R_{k}\lvert^{2})\rvert R_{k}\lvert^{2}dx=\int _{\Omega_{1}}E_{1}(\rvert Q_{\omega_{k}}\lvert^{2})\rvert Q_{\omega_{k}}\lvert^{2}dx \leq Ce^{-2\theta_{0}t},
\end{align*}
\begin{align*}
\int _{\Omega_{1}}E_{1}(Re(R_{k}\overline{R}_{k}))Re(R_{k}\overline{R}_{k})dx\leq Ce^{-2\theta_{0}t},
\end{align*}
\begin{align*}
\int _{\Omega_{1}}E_{1}(\rvert \varepsilon\lvert^{2})\rvert \varepsilon\lvert^{2}dx\leq C\|\varepsilon\|_{L^{4}}^{4}\leq C\|\varepsilon\|_{H^{1}}^{4}\leq Ce^{-2\theta_{0}t},
\end{align*}
\begin{align*}
\int _{\Omega_{1}}E_{n}(\rvert R_{k}\lvert^{2})^{2}dx=\int _{\Omega_{1}}E_{n}(\rvert Q_{\omega_{k}}\lvert^{2})^{2}dx \leq Ce^{-2\theta_{0}t}.
\end{align*}
Hence we have
\begin{align}
\int _{\Omega_{1}}E_{1}(\rvert \varphi_{n}\lvert^{2})\rvert \varphi_{n}\lvert^{2}+\rvert \nabla z_{n} \lvert^{2}dx \leq Ce^{-2\theta_{0}t}
\end{align}
Combining with (5.49)-(5.52), the support properties of $Y$ and Sobolev imbedding we obtain
\begin{align}
\big\rvert \frac{d}{dt}Im\int\partial_{x_{1}}\varphi_{n}\bar \varphi_{n} Ydx\big\rvert\leq \frac{C}{L}\int_{\Omega_{1}}(\rvert \nabla \varphi_{n}\lvert^{2}+\rvert\varphi_{n}\lvert^{2}+\rvert\varphi_{n}\lvert^{p+1})dx.%\leq \frac{C}{L}\int(|\nabla \varphi_{n}|^{2}+|\varphi_{n}|^{2}+|\varphi_{n}|^{p+1})dx.
\end{align}
Now by the Sobolev inequality applied to $\varphi_{n}(x)h(x_{1}-\sigma_{k}t)$, where $h=h(x_{1})$ is a $C^{1}$ function such that $h(x_{1})=1$ for $\rvert x_{1}\lvert<L$ and $h(x_{1})=0$ for $\rvert x_{1}\lvert>L+1$, we have
\begin{eqnarray}
\int_{\Omega_{1}}\rvert\varphi_{n}\lvert^{p+1}dx\leq C(\int_{\tilde{\Omega}_{1}} \rvert \varphi_{n}\lvert^{2}+\rvert \nabla \varphi_{n}\lvert^{2}dx)^{\frac{p+1}{2}},
\end{eqnarray}
where
\begin{eqnarray*}
\tilde{\Omega}_{1}(t)=[-(L+1)+\sigma_{k}t,(L+1)+\sigma_{k}t]\times \mathbb{R}^{1}.
\end{eqnarray*}
From (5.53) and (5.54), we obtain
\begin{align}
\rvert \frac{d}{dt}Im\int\partial_{x_{1}}\varphi_{n}\bar \varphi_{n} Ydx\lvert\leq \frac{C}{L}\int_{\tilde{\Omega}_{1}}(\rvert\nabla \varphi_{n}\lvert^{2}+\rvert\varphi_{n}\lvert^{2}+(\rvert\nabla \varphi_{n}\lvert^{2}+\rvert\varphi_{n}\lvert^{2})^{\frac{p+1}{2}})dx.\notag\\
\end{align}
By (5.41), we have
\begin{align}
\frac{1}{2}\frac{d}{dt}Im\int\partial_{x_{2}}\varphi_{n}\bar \varphi_{n}Y)dx&=Re\int\partial_{x_{2}}\varphi_{n}\partial_{x_{1}}\bar \varphi_{n}Y'dx
+\frac{1}{2}\int\partial_{x_{1}}z_{n}\partial_{x_{2}}z_{n}Y'dx\notag\\
&-\frac{1}{2}\int\partial_{x_{2}}z_{n}\cdot \rvert \varphi_{n}\lvert^{2}Y'dx-\frac{\sigma_{k}}{2L}Im\int \partial_{x_{2}}\varphi_{n} \bar \varphi_{n} Y'dx.
\end{align}
Similar arguments to those as before, we have
\begin{align}
\frac{1}{2}\frac{d}{dt}Im\int\partial_{x_{2}}\varphi_{n}\bar \varphi_{n}Ydx \leq \frac{C}{L}\int_{\Omega_{1}}(\rvert \nabla \varphi_{n}\lvert^{2}+\rvert\varphi_{n}\lvert^{2})dx.
\end{align}
\par  Next,  by $\varphi_{n}(t)=R(t)+(\varphi_{n}(t)-R(t))$, we have
\begin{align}
\int_{\widetilde{\Omega}_{1}}(\rvert\nabla \varphi_{n}(t)\lvert^{2}+\rvert \varphi_{n}(t)\lvert^{2})dx&\leq 2\int_{\widetilde{\Omega}_{1}}(\rvert\nabla R(t)\lvert^{2}+\rvert R(t)\lvert^{2})dx\notag\\
&+2\|\varphi_{n}(t)-R(t)\|_{H^{1}}^{2}.
\end{align}
By Lemma 2.4, $Q_{\omega}$ has exponential decay property
\begin{eqnarray*}
\rvert\nabla Q_{\omega}(x)\lvert+\rvert Q_{\omega}(x)\lvert\leq Ce^{-\frac{\sqrt{\omega}}{2}\rvert x\lvert}.
\end{eqnarray*}
Thus by the definition of $\theta_{0}$ and $\sigma_{k}$, we can make the following conclusion
\begin{eqnarray}
\int_{\widetilde{\Omega}_{1}}(\rvert\nabla R(t)\lvert^{2}+\rvert R(t)\lvert^{2})dx\leq Ce^{-8\sqrt{\theta_{0}}(\sqrt{\theta_{0}}t-L)}\leq Ce^{-4\theta_{0}t}
\end{eqnarray}
by taking $T_{0}$ and $L$ such that $\sqrt{\theta_{0}}T_{0}\geq 2L$. Therefore, from (5.21), (5.48), (5.55)-(5.59) and the definition of $I_{k}(t)$ and $M_{k}(t)$, and taking $A_{0}e^{-\theta_{0}T_{0}}$ small enough, we have
\begin{equation}
\rvert\frac{d}{dt}I_{k}(t)\lvert+\rvert\frac{d}{dt}M_{k}(t)\lvert\leq \frac{CA_{0}^{2}}{L}e^{-2\theta_{0}t}.
\end{equation}
Note that for $I_{1}(t)$ and $M_{1}(t)$ we have also used the conservations of mass and momentum. Now by integrating (5.60) between $t$ and $T_{n}$, we obtain
\begin{eqnarray*}
\rvert I_{k}(T_{n})-I_{k}(t)\lvert+\rvert M_{k}(T_{n})-M_{k}(t)\lvert\leq \frac{CA_{0}^{2}}{L}e^{-2\theta_{0}t}.
\end{eqnarray*}
\par This completes the proof of Lemma 5.11.
\par \noindent\textbf{Lemma 5.12.} There exists $C>0$ such that for any $t\in[t^{*},\;T_{n}]$,
\begin{equation}
\lvert\omega_{k}(t)-\omega^{0}_{k}\rvert\leq C\|\varepsilon(t)\|^{2}_{L^{2}}+C(\frac{A_{0}^{2}}{L}+1)e^{-2\theta_{0}t}.
\end{equation}
\noindent\textbf{Proof.} From (5.31) and (5.45), we have
\begin{eqnarray*}
I_{k}(t)=\int\lvert \tilde{R}(t)\rvert^{2}y_{k}(t)dx+2Re \int\tilde{R}(t)\overline{\varepsilon}(t)y_{k}(t)dx+\int\lvert\varepsilon(t)\rvert^{2}y_{k}(t)dx.
\end{eqnarray*}
By the exponential decay of each $Q_{\omega_{k}(t)}$, the orthogonality $\int\tilde{R}_{k}(t)\bar{\varepsilon}(t)dx=0$ and the property of support of $y_{k}$, we have
\begin{eqnarray*}
I_{k}(t)=\int\lvert\varphi_{n}(t)\rvert^{2}y_{k}(t)dx=\int Q^{2}_{\omega_{k}(t)}dx+\int \lvert\varepsilon(t)\rvert^{2}y_{k}(t)dx+O(e^{-2\theta_{0}t}).
\end{eqnarray*}
From the result of Lemma 5.11, we have
\begin{eqnarray*}
\lvert I_{k}(t)-I_{k}(T_{n})\rvert\leq \frac{CA_{0}^{2}}{L}e^{-2\theta_{0}t}.
\end{eqnarray*}
Thus, by $\omega_{k}(T_{n})=\omega_{k}^{0}$ and $\varepsilon(T_{n})\equiv 0$, we obtain
\begin{eqnarray}
\lvert\int Q^{2}_{\omega_{k}(t)}dx-\int Q^{2}_{\omega_{k}^{0}}dx\rvert\leq C\|\varepsilon(t)\|^{2}_{L^{2}}+C(\frac{A_{0}^{2}}{L}+1)e^{-2\theta_{0}t}.
\end{eqnarray}
Recall that $\frac{d}{d\omega}\int Q^{2}_{\omega}dx\lvert_{\omega=\omega_{k}^{0}}>0$, then we assume $\omega_{k}(t)$ is close to $\omega_{k}^{0}$. Thus
\begin{align*}
(\omega_{k}(t)-\omega^{0}_{k})(\frac{d}{d\omega}\int Q^{2}_{\omega}dx\lvert_{\omega=\omega^{0}_{k}})
&=\int Q^{2}_{\omega_{k}(t)}dx-\int Q^{2}_{\omega_{k}^{0}}dx\notag\\
&-\beta(\omega_{k}(t)-\omega^{0}_{k})(\omega_{k}(t)-\omega^{0}_{k})^{2}
\end{align*}
with $\beta(\epsilon)\rightarrow 0$, as $\epsilon\rightarrow 0$, which implies that for some constant $C=C(\omega_{k}^{0})$.
\begin{eqnarray}
\lvert\omega_{k}(t)-\omega^{0}_{k}\rvert\leq C \lvert\int Q^{2}_{\omega_{k}(t)}dx-\int Q^{2}_{\omega_{k}^{0}}dx\rvert.
\end{eqnarray}
Therefore by (5.62) and (5.63), we have
\begin{eqnarray*}
\lvert\omega_{k}(t)-\omega^{0}_{k}\rvert\leq C\|\varepsilon(t)\|^{2}_{L^{2}(\mathbb{R}^{2})}+C(\frac{A_{0}^{2}}{L}+1)e^{-2\theta_{0}t}.
\end{eqnarray*}
\par  This proves Lemma 5.12.
\par \noindent\textbf{Lemma 5.13.} Let $1 < p < 3$ and $\omega_{k}^{0}\in (0,\;\omega_{J})$. Then there exists $\lambda > 0$ such that for any real-valued $v\in H^{1}$ satisfying $Re(Q_{\omega_{k}^{0}},\;v)=Im(Q_{\omega_{k}^{0}},\;v)=0$ and $Re(\nabla Q_{\omega_{k}^{0}},\;v)=0$, one has that
\begin{eqnarray}
(H_{\omega_{k}^{0}}v,v)\geq \lambda\|v\|_{H^{1}}^{2}.
\end{eqnarray}
\noindent\textbf{Proof.}  By (4.19) and Theorem 4.6, we have that $D''(\omega^{0}_{k}) > 0$. From Theorem 3.3 and Corollary 3.31 in \cite{GSS1987}, we get this result.
\par \noindent\textbf{Lemma 5.14.} Let $1 < p < 3$. For $\omega_{k}^{0}\in (0,\;\omega_{J})$ and $\omega_{k}(t)$ close to $\omega_{k}^{0}$, we have
\begin{eqnarray*}
\lvert\Gamma_{\omega_{k}^{0}}(Q_{\omega_{k}(t)})-\Gamma_{\omega_{k}^{0}}(Q_{\omega_{k}^{0}})\rvert\leq C\lvert\omega_{k}(t)-\omega_{k}^{0}\rvert^{2},
\end{eqnarray*}
where $\Gamma_{\omega_{k}^{0}}(z)=E(z)+\omega_{k}^{0}M(z)$.
\par \noindent\textbf{Proof.} By (2.2) and (2.3), we have
\begin{align}
\Gamma_{\omega_{k}^{0}}(Q_{\omega_{k}(t)})=E(Q_{\omega_{k}(t)})+\omega_{k}^{0}\int \lvert Q_{\omega_{k}(t)}\rvert^{2}dx.
\end{align}
By Taylar expansion of $\Gamma_{\omega_{k}^{0}}(Q_{\omega_{k}(t)})$, (5.65), Theorem 4.2 and Theorem 4.6, we have
\begin{align}
\Gamma_{\omega_{k}^{0}}(Q_{\omega_{k}(t)})
%={} &\Gamma_{\omega_{k}^{0}}(Q_{\omega_{k}^{0}})+\frac{d}{d\omega}\Gamma_{\omega_{k}^{0}}(Q_{\omega})\lvert_{\omega=\omega_{k}^{0}}\cdot (\omega_{k}(t)-\omega_{k}^{0})\notag\\
%+&\lvert\omega_{k}(t)-\omega_{k}^{0}\rvert^{2}\beta(\lvert\omega_{k}(t)-\omega_{k}^{0}\rvert) \notag\\
%={} &\Gamma_{\omega_{k}^{0}}(Q_{\omega_{k}^{0}})+\lvert\omega_{k}(t)-\omega_{k}^{0}\rvert^{2}\beta(\lvert\omega_{k}(t)-\omega_{k}^{0}\rvert) \notag\\
%+{} &[\frac{d}{d\omega}E(Q_{\omega_{k}(t)})\lvert_{\omega_{k}(t)=\omega_{k}^{0}}+\omega_{k}^{0}\frac{d}{d\omega}\int Q^{2}_{\omega}dx\lvert_{\omega=\omega_{k}^{0}}]\cdot (\omega_{k}(t)-\omega_{k}^{0})\notag\\
={} &\Gamma_{\omega_{k}^{0}}(Q_{\omega_{k}^{0}})-(\omega_{k}(t)-\omega_{k}^{0})^{2}\frac{d}{d\omega}\int Q^{2}_{\omega}dx\lvert_{\omega=\omega_{k}^{0}}\notag\\
+&\lvert\omega_{k}(t)-\omega_{k}^{0}\rvert^{2}\beta(\lvert\omega_{k}(t)-\omega_{k}^{0}\rvert).
\end{align}
By (5.66), Theorem 4.6 and $\omega_{k}(t)$ close to $\omega^{0}_{k}$, there exists $C=C(\omega^{0}_{k}) > 0$ such that
\begin{eqnarray*}
\lvert\Gamma_{\omega_{k}^{0}}(Q_{\omega_{k}(t)})-\Gamma_{\omega_{k}^{0}}(Q_{\omega_{k}^{0}})\rvert\leq C\lvert\omega_{k}(t)-\omega_{k}^{0}\rvert^{2}.
\end{eqnarray*}
\par  This completes the proof of Lemma 5.13.
\par  Now we set
\begin{equation}
J(t)=\sum_{k=1}^{K}[(\omega^{0}_{k}+\frac{1}{4}\lvert v_{k}\rvert^{2})I_{k}(t)-v_{k}M_{k}(t)]
\end{equation}
and
\begin{equation}
G(t)=E(\varphi_{n}(t))+J(t).
\end{equation}
\par  From (5.43) to (5.45), Lemma 2.4 and Lemma 5.13, Lemma 5.14 directly deduces the following Lemma.
\par \noindent\textbf{Lemma 5.15.}
For all $t\in[t^{*},T_{n}]$, we have
\begin{align}
G(t)=&\sum_{k=1}^{K}[E(Q_{\omega_{k}^{0}})+\omega_{k}^{0}\int Q^{2}_{\omega_{k}^{0}}dx]+P(\varepsilon(t),\varepsilon(t))+\sum_{k=1}^{K}O(\lvert\omega_{k}(t)-\omega_{k}^{0}\rvert^{2})\notag\\
+&\|\varepsilon(t)\|_{H^{1}}^{2}\beta (\|\varepsilon(t)\|_{H^{1}})+O(e^{-2\theta_{0}t})
\end{align}
with $\beta(\epsilon)\rightarrow 0$, as $\epsilon\rightarrow 0$, where
\begin{align}
P(\varepsilon,\varepsilon)
={} &\int\lvert\nabla\varepsilon\rvert^{2}dx-\sum_{k=1}^{K}(\int\lvert\widetilde{R}_{k}\rvert^{p-1}\lvert\varepsilon\rvert^{2}+(p-1)\lvert\widetilde{R}_{k}\rvert^{p-3}(Re( \overline{\widetilde{R}}_{k}\varepsilon))^{2}dx)\notag\\
+&\sum_{k=1}^{K}((\omega_{k}(t)+\frac{1}{4}\lvert v_{k}\rvert^{2})\int \lvert\varepsilon\rvert^{2} y_{k}(t)dx-v_{k}\cdot Im\int \nabla \varepsilon\cdot\overline{\varepsilon}y_{k}(t)dx)\notag\\
-&\sum_{k=1}^{K}\frac{1}{2}\int (E_{1}(\lvert\widetilde{R}_{k}\rvert^{2})\lvert\varepsilon\rvert^{2}+ E_{1}(\lvert\varepsilon\rvert^{2})\lvert\widetilde{R}_{k}\rvert^{2}+4 E_{1}( Re(\overline{\widetilde{R}}_{k}\varepsilon))\widetilde{R}_{k}\varepsilon)dx.
\end{align}
\noindent\textbf{Proof.} For $\omega_{k}(t),\;\;\omega^{0}_{k}\in (0,\;\omega_{J})$ and $\omega_{k}(t)$ close to $\omega^{0}_{k}$, from Lemma 5.14, we have that
\begin{equation}
\rvert E(Q_{\omega^{0}_{k}})+\omega^{0}_{k}\int Q_{\omega^{0}_{k}}^{2}dx-E(Q_{\omega_{k}(t)})-\omega^{0}_{k}\int Q_{\omega_{k}(t)}^{2}dx\lvert \leq C\rvert \omega_{k}(t)-\omega^{0}_{k}\lvert^{2}.
\end{equation}
%\par We refer to Weinstein[8], Section 2, Eq.(2.5) for this property.
Now, by the definition of $y_{k}$, (5.67) and (5.68), we have $\sum_{k=1}^{K} y_{k}=1$. Thus
\begin{align}
G(t)=&\sum_{k=1}^{K}\int(\rvert\nabla \varphi_{n}\lvert^{2}-\frac{1}{2}E_{1}(\lvert\varphi_{n}\rvert^{2})\rvert\varphi_{n}\lvert^{2}-\frac{2}{p+1}\rvert\varphi_{n}\lvert^{p+1}\notag\\
+&(\omega^{0}_{k}+\frac{1}{4}\rvert v_{k}\lvert^{2})\rvert\varphi_{n}\lvert^{2}-v_{k}Im(\nabla\varphi_{n}\overline{\varphi}_{n}))y_{k}dx.
\end{align}
Expanding $\varphi_{n}(t)=\widetilde{R}(t)+\varepsilon(t)$ in the expression of $E(\varphi_{n}(t))$. By the calculations, we have that
\begin{align}
E(\varphi_{n})
{} &=E(\widetilde{R})-2Re\int(\Delta \overline{\widetilde{R}}+\rvert\widetilde{R}\lvert^{p-1} \overline{\widetilde{R}}+E_{1}(\rvert\widetilde{R}\lvert^{2}) \overline{\widetilde{R}})\varepsilon dx\notag\\
&-\int\rvert\widetilde{R}\lvert^{p-1}\rvert\varepsilon\lvert^{2}+(p-1)\rvert\widetilde{R}\lvert^{p-3}(Re( \overline{\widetilde{R}}\varepsilon))^{2}dx\notag\\
&-\frac{1}{2}\int E_{1}(\lvert\widetilde{R}\rvert^{2})\lvert\varepsilon\rvert^{2}+ E_{1}(\lvert\varepsilon\rvert^{2})\lvert\widetilde{R}\rvert^{2}+4 E_{1}( Re(\overline{\widetilde{R}}\varepsilon))\overline{\widetilde{R}}\varepsilon))\notag\\
&-E_{1}(Re(\overline{\widetilde{R}}\varepsilon ))\rvert\widetilde{R}\lvert^{2} +\|\varepsilon\|_{H^{1}}^{2}\beta(\|\varepsilon\|_{H^{1}}).
\end{align}
Note that the $\widetilde{R}_{k}(t)$ and $E_{1}(\rvert\widetilde{R}_{k}(t)\lvert^{2})$ are exponentially decaying, we have that
\begin{align}
E(\varphi_{n})
{} &=\sum^{K}_{k=1}(E(\widetilde{R}_{k})-2Re\int(\Delta \overline{\widetilde{R}}_{k}+\rvert\widetilde{R}_{k}(t)\lvert^{p-1}\overline{\widetilde{R}}_{k}+E_{1}(\rvert\widetilde{R}_{k}\lvert^{2}) \overline{\widetilde{R}}_{k})\varepsilon  dx)\notag\\
&-\sum^{K}_{k=1}\int\rvert\widetilde{R}_{k}\lvert^{p-1}\rvert\varepsilon\lvert^{2}+(p-1)\rvert\widetilde{R}_{k}\lvert^{p-3}(Re(\overline{\widetilde{R}}_{k}\varepsilon))^{2}dx\notag\\
-&\sum_{k=1}^{K}\frac{1}{2}\int (E_{1}(\lvert\widetilde{R}_{k}\rvert^{2})\lvert\varepsilon\rvert^{2}+ E_{1}(\lvert\varepsilon\rvert^{2})\lvert\widetilde{R}_{k}\rvert^{2}+4 E_{1}( Re(\overline{\widetilde{R}}_{k}\varepsilon))\overline{\widetilde{R}}_{k}\varepsilon)dx\notag\\
&+\int\rvert\nabla\varepsilon\lvert^{2}dx+\|\varepsilon\|_{H^{1}}^{2}\beta(\|\varepsilon\|_{H^{1}})+O(e^{-2\theta_{0}t}).
\end{align}
Now we turn to $J(t)$. Expanding $\varphi_{n}(t)=\widetilde{R}(t)+\varepsilon(t)$ in the expression of $I_{k}(t)$
\begin{align*}
I_{k}(t)=\int\rvert\widetilde{R}(t)\lvert^{2}y_{k}(t)dx+\int\rvert\varepsilon(t)\lvert^{2}y_{k}(t)dx+2Re\int \overline{\widetilde{R}}(t)\varepsilon(t)y_{k}(t)dx.
\end{align*}
By the properties of $y_{k}$, the properties of $\widetilde{R}(t)$ and the orthogonality conditions on $\varepsilon(t)$, we get that
\begin{align*}
I_{k}(t)=\int\rvert\widetilde{R}_{k}(t)\lvert^{2}dx+\int\rvert\varepsilon(t)\lvert^{2}y_{k}(t)dx+O(e^{-2\theta_{0}t}).
\end{align*}
Similarly, for $M_{k}(t)$, we have
\begin{align*}
M_{k}(t)=Im\int\nabla \widetilde{R}_{k}\overline{\widetilde{R}}_{k}dx-2Im\int\nabla \overline{\widetilde{R}}_{k}\varepsilon dx+Im\int\nabla\varepsilon\overline{\varepsilon} y_{k}(t)dx+O(e^{-2\theta_{0}t}).
\end{align*}
It follows that
\begin{align}
J(t)&=\sum_{k=1}^{K}((\omega^{0}_{k}+\frac{1}{4}\rvert v_{k}\lvert^{2})(\int\rvert\widetilde{R}_{k}\lvert^{2}dx+2Re\int \overline{\widetilde{R}}_{k}\varepsilon dx+\int\rvert\varepsilon\lvert^{2}y_{k}(t)dx))\notag\\
&-\sum_{k=1}^{K}(v_{k}(Im\int\nabla\widetilde{R}_{k}\overline{\widetilde{R}}_{k}dx-2Im\int\nabla\overline{\widetilde{R}}_{k}\varepsilon  dx+Im\int\nabla\varepsilon\overline{\varepsilon}y_{k}(t)dx))\notag\\
&+O(e^{-2\theta_{0}t}).
\end{align}
By the equation of $\widetilde{R}_{k}(t)$, and the orthogonality conditions on $\varepsilon(t)$, we have
\begin{align*}
-2Re\int(\Delta \overline{\widetilde{R}}_{k}+\rvert\widetilde{R}_{k}\lvert^{p-1} \overline{\widetilde{R}}_{k}+E_{1}(\rvert\widetilde{R}_{k}\lvert^{2}) \overline{\widetilde{R}}_{k})\varepsilon dx+&2(\omega^{0}_{k}+\frac{1}{4}\rvert v_{k}\lvert^{2})Re\int\overline{\widetilde{R}}_{k}\varepsilon dx\notag\\
+&2v_{k}Im\int\nabla\overline{\widetilde{R}}_{k}\varepsilon dx=0,
\end{align*}
which means that the terms of order 1 in $\varepsilon(t)$ all disappear when we sum $E(\varphi_{n}(t))$ and $J(t)$.
Therefore, with the definition of $P(\varepsilon(t),\varepsilon(t))$, we obtain (5.69).
\par  This completes the proof of Lemma 5.15.
\par \noindent\textbf{Lemma 5.16.}(Lemma 4.11 in \cite{MMT2006}) There exists $\lambda>0$ such that for all $t\in[t^{*},T_{n}]$,
\begin{eqnarray}
P(\varepsilon(t),\varepsilon(t))\geq \lambda\|\varepsilon(t)\|_{H^{1}}^{2}.
\end{eqnarray}
\par  Combining with Lemma 5.10, Lemma 5.11, Lemma 5.15 and Lemma 5.16, we can deduce the following lemma according to Martel and Merle's way \cite{MM2006}.
\par \noindent\textbf{Lemma 5.17.} (Lemma 5 in \cite{MM2006}) For any $t\in[t^{*},T_{n}]$
\begin{align}
\|\varepsilon(t)\|_{H^{1}}^{2}+\lvert\omega_{k}(t)-\omega_{k}^{0}\rvert+\lvert x_{k}(t)\rvert^{2}+\lvert\gamma_{k}(t)-\gamma_{k}^{0}\rvert^{2}\leq C(\frac{A_{0}^{2}}{L}+1)e^{-2\theta_{0}t}.
\end{align}
\par \noindent\textbf{Lemma 5.18.}  For any $t\in [t^{*},T_{n}]$, there exists $C > 0$ such that
\begin{equation}
\|R(t)-\tilde{R}(t)\|_{H^{1}(\mathbb{R}^{2})}\leq C\sum_{k=1}^{K}(\lvert\omega_{k}(t)-\omega_{k}^{0}\rvert+\lvert x_{k}(t)\rvert+\lvert\gamma_{k}(t)-\gamma_{k}^{0}\rvert).
\end{equation}
\noindent\textbf{Proof.} By (5.5), (5.9), (5.32), (5.33) and (5.34), we have
\begin{align}
\widetilde{R}_{k}(t)=&R_{k}(t)+(\omega_{k}(t)-\omega^{0}_{k})\frac{dQ_{\omega_{k}(t)}}{d\omega}(\cdot-\widetilde{x}_{k}(t))e^{i(\frac{1}{2}v_{k}x+\delta_{k}(t))}\lvert_{\omega_{k}(t)=\omega^{0}_{k},x_{k}(t)=0,\gamma_{k}(t)=\gamma^{0}_{k}}\notag\\
-&x_{k}(t)\nabla\widetilde{R}_{k}(t)\lvert_{\omega_{k}(t)=\omega^{0}_{k},x_{k}(t)=0,\gamma_{k}(t)=\gamma^{0}_{k}}+i(\gamma_{k}(t)-\gamma^{0}_{k})\widetilde{R}_{k}(t)\lvert_{\omega_{k}(t)=\omega^{0}_{k},x_{k}(t)=0,\gamma_{k}(t)=\gamma^{0}_{k}}\notag\\
+&O((\omega_{k}(t)-\omega^{0}_{k})^{2})+O(x^{2}_{k}(t))+O((\gamma_{k}(t)-\gamma^{0}_{k})^{2}).
\end{align}
By (5.79), Lemma 5.9 and Lemma 5.17 deduce that
\begin{align*}
\|R(t)-\tilde{R}(t)\|_{H^{1}}\leq C\sum_{k=1}^{K}(\lvert\omega_{k}(t)-\omega_{k}^{0}\rvert+\lvert x_{k}(t)\rvert+\lvert\gamma_{k}(t)-\gamma_{k}^{0}\rvert).
\end{align*}
\par  This proves Lemma 5.18.

\noindent\textbf{Proof of Proposition 5.7.} From Lemma 5.18, we get for all $t\in[t^{*},T_{n}]$
\begin{align}
\|R(t)-\tilde{R}(t)\|_{H^{1}}^{2}\leq& C\sum_{k=1}^{K}(\lvert\omega_{k}(t)-\omega_{k}^{0}\rvert^{2}+\lvert\gamma_{k}(t)-\gamma_{k}^{0}\rvert^{2}+\lvert x_{k}(t)\rvert^{2})\notag\\
\leq& C(\frac{A_{0}^{2}}{L}+1)e^{-2\theta_{0}t},
\end{align}
By Lemma 5.17 and (5.80), we have
\begin{eqnarray*}
\|\varphi_{n}(t)-R(t)\|_{H^{1}}^{2}\leq 2\|\varepsilon(t)\|_{H^{1}}^{2}+2\|\tilde{R}(t)-R(t)\|_{H^{1}}^{2}\leq C(\frac{A_{0}^{2}}{L}+1)e^{-2\theta_{0}t},
\end{eqnarray*}
where $C>0$ does not depend on $A_{0}$. Now we choose $A_{0}^{2}>8C$, $L=A_{0}^{2}$,
and $T_{0}$ large enough. It follows that
\begin{eqnarray*}
\|\varphi_{n}(t)-R(t)\|_{H^{1}}^{2}\leq 2Ce^{-2\theta_{0}t}\leq \frac{A_{0}^{2}}{4}e^{-2\theta_{0}t}.
\end{eqnarray*}
Therefore, the conclusion is that for any $t\in[t^{*},T_{n}]$, $\|\varphi_{n}(t)-R(t)\|_{H^{1}}\leq \frac{A_{0}}{2}e^{-\theta_{0}t}$.
\par  This completes the proof of Proposition 5.7.
\par \noindent\textbf{Corollary 5.19.} For multi-solitons  $\varphi(t, x)$ of (1.1) in Theorem 5.1, we have that $\varphi(t, x)$ satisfying $\int \lvert\varphi(t, x)\rvert^{2}dx < 2d_{J}$ with $t\in \mathbb{R}$.
\\
\noindent\textbf{Proof.} From Claim 5.5,
\begin{align*}
\|\psi_{0}\|_{L^{2}(\mathbb{R}^{2})}\leq \liminf_{n\rightarrow\infty}\|\varphi_{n}(T_{0})\|_{L^{2}}<\sqrt{2d_{J}}.
\end{align*}
By Theorem 5.1,
\begin{align*}
\|\varphi(t)\|_{L^{2}}=\|\psi_{0}\|_{L^{2}},\;\;t\in \mathbb{R}.
\end{align*}
It follows that $\int \lvert\varphi(t)\rvert^{2}dx < 2d_{J} $ for  $t\in \mathbb{R}$.
\par  This proves Corollary 5.19.

\par \noindent\textbf{Acknowledgment.}
\par  This research is supported by the National Natural Science Foundation of China 11871138.

\end{document}